\documentclass{amsart}
\usepackage{amssymb,amsmath,epsfig}
\usepackage[all]{xy}

\theoremstyle{plain}
\newtheorem{thm}[equation]{Theorem}
\newtheorem{prop}[equation]{Proposition}
\newtheorem{lemma}[equation]{Lemma}
\newtheorem{cor}[equation]{Corollary}
\newtheorem{remark}[equation]{Remark}
\newtheorem{example}[equation]{Example}
\newtheorem{defn}[equation]{Definition}

\numberwithin{equation}{section}

\newcommand{\kk}{{\mathbf k}}
\newcommand{\kku}{\kk[u]/u^n} 
\newcommand{\kkux}{(\kku)^{\times}}

\newcommand{\kkuek}{ \kk[u]/u^{e_K l} }
\newcommand{\kkuekx}{(\kkuek)^{\times}} 
 
\newcommand{\HH}{H}
\newcommand{\gpi}{\left(\frac{g \pi}{\pi}\right)}

\newcommand{\rhobar}{\overline{\rho}}
\newcommand{\M}{{\mathcal M}}
\newcommand{\MM}{\M_1}
\newcommand{\N}{{\mathcal N}}
\newcommand{\NN}{\N_1}
\newcommand{\p}{\phi_1}
\newcommand{\unif}{(-l)^{\frac{1}{l^2-1}}}
\newcommand{\e}{ {\rm \bf e} }
\newcommand{\etil}{ \tilde{\e} }

\newcommand{\ee}{ {\rm \bf e}'}
\newcommand{\f}{ { {\rm \bf f} \, } }
\newcommand{\ff}{ { {\rm \bf f} \,'} }
\newcommand{\Id}{\text{Id}}
\newcommand{\Ql}{{\mathbb Q}_l}
\newcommand{\Z}{{\mathbb Z}}

\newcommand{\Zl}{\Z_l}
\newcommand{\Zlx}{\Z_l^{\times}}
\newcommand{\Zll}{\Z_{l^2}}
\newcommand{\Q}{{\mathbb Q}}

\newcommand{\Qlbar}{\overline{\Q}_l}
\newcommand{\Qll}{{\mathbb Q}_{l^2}}

\newcommand{\F}{{\mathbb F}}
\newcommand{\Fl}{{\mathbb F}_l}
\newcommand{\Fll}{{\mathbb F}_{l^2}}
\newcommand{\Flbar}{\overline{\F}_l}
\newcommand{\Flx}{\Fl^{\times}}
\newcommand{\Fllx}{\Fll^{\times}}

\newcommand{\Frob}{{\rm Frob}}
\newcommand{\Fr}{{\rm Frob_l}}
\newcommand{\OO}{\mathcal O}
\newcommand{\G}{\mathcal G}
\newcommand{\la}{\langle}
\newcommand{\ra}{\rangle}
\newcommand{\Gal}{{\rm Gal}}
\newcommand{\Hom}{{\rm Hom}}
\newcommand{\Ext}{{\rm Ext}}
\newcommand{\Trace}{{\rm Trace}}

\newcommand{\coeffFll}{\Fll[u]/u^{l(l^2-1)}}
\newcommand{\z}{[\zeta]}
\newcommand{\fr}{[\varphi]}
\newcommand{\GL}{{\rm GL}_2}

\newcommand{\om}{\tilde{\omega}}
\newcommand{\omt}{\tilde{\omega}_2}

\newcommand{\Rdvz}{R(\rhobar,2,\tau)}
\newcommand{\s}{\mathcal S}

\DeclareMathOperator{\Spec}{{\rm Spec \ }}

\begin{document}
\title{Modularity of Some Potentially Barsotti-Tate Galois
Representations}
\author{David Savitt}
\date{May, 2002}
\maketitle

\pagestyle{myheadings}
\markboth{DAVID SAVITT}{MODULARITY OF SOME POTENTIALLY BARSOTTI-TATE...}

\section{Notation, terminology, and results}

\noindent Throughout this article, we let $l$ be an odd prime, and we fix 
an algebraic closure~$\Qlbar$ of~$\Ql$ with residue field $\Flbar$.  The 
fields $K$, $L$, and $E$ will 
always be finite extensions of~$\Ql$ inside~$\Qlbar$.  We denote 
by~$G_{E}$ the 
Galois group~$\Gal(\Qlbar/E)$, by~$W_E$ the Weil group of~$E$, and 
by~$I_E$ the inertia group of $E$.  The group~$I_{\Ql}$ will be 
abbreviated~$I_l$.  The character $\omega_n : G_{\Ql} 
\rightarrow \F_{l^n} \subset \Flbar$ is defined via
$$
\omega_n : u \mapsto \frac{ 
u\left(\left(-l\right)^{\frac{1}{l^n-1}}\right) } 
{\left(-l\right)^{\frac{1}{l^n-1}}} \,,
$$
and its Teichm\"uller lift will be denoted~$\om_n$.  In particular, 
$\omega = 
\omega_1$ is the mod-$l$ reduction of the cyclotomic character~$\epsilon$.  
Recall that if $\rho 
: G_{\Ql} \rightarrow \GL(K)$  or $W_{\Ql} \rightarrow \GL(K)$
is continuous and 
tamely 
ramified, then $\rho \,|_{I_l} \otimes_K \overline{K}$ is isomorphic 
either to $\om^i \oplus 
\om^j$ or to $\omt^m \oplus \omt^{lm}$, depending on the absolute 
reducibility or irreducibility of~$\rho$.  

If an $l$-adic representation~$\rho$ of~$G_{\Ql}$ is potentially 
semistable (in the sense of Fontaine), then one associates to~$\rho$ a 
Weil-Deligne representation $WD(\rho)$ over $\Qlbar$, 
for example as in Section B.1 
of \cite{CDT}.  Then $\rho$ becomes semistable over $E$ if and only if 
$WD(\rho) \,|_{I_E}$ is trivial.  The \textit{Galois type}~$\tau(\rho)$ 
associated to such $\rho$ is defined to be the the isomorphism class of 
the representation $WD(\rho)\,|_{I_l}$ of~$I_l$.  

Following \cite{BCDT} and using the notation of \cite{BreuilMezard}, we 
define a 
collection of deformation rings.  Let $\overline{\rho} 
: G_{\Ql}
\rightarrow {\rm GL}_2(\kk)$ be a representation over a finite
field~$\kk$ of characteristic~$l$, and assume that the only matrices which 
commute with the image of~$\overline{\rho}$ are scalar matrices.  Fix
a positive integer~$k$, and let~$\tau$ be a Galois
type.
We are
interested in lifts $\rho : G_{\Ql} \rightarrow {\rm GL}_2 (\Qlbar)$ 
of~$\overline{\rho}$
with the following properties:
\begin{enumerate} 
\item $\rho$ is potentially semi-stable with Hodge-Tate
weights $(0,k-1)$,
\item $\tau(\rho)$ is isomorphic to $\tau$, and
\item $\det(\rho) = \epsilon^{k-1} \chi$,
 where $\chi$ is a character
of finite order prime to $l$.
\end{enumerate}
Let $R^{univ}_{\OO}$ denote the universal deformation 
ring parametrizing deformations of $\rhobar$ over complete 
local Noetherian $\OO$-algebras, 
where $\OO$ is the integer ring of a finite extension of $\Ql$ inside 
$\Qlbar$ which contains both 
the 
Witt 
vectors 
$W(\kk)$ and a field of rationality
of $\tau$.
Let $\rho^{univ}$ be the universal deformation.  
We say that a prime ${\mathfrak p}$ 
of $R^{univ}_{\OO}$ has type $(k,\tau)$ if there is a field $K \supset
\OO$ 
and a map of $\OO$-algebras
$$ f_{{\mathfrak p}} : R^{univ}_{\OO} \rightarrow K \ \ \text{with} \ \  
{\mathfrak p} = \ker( f_{{\mathfrak p}} )$$ such that the
pushforward of $\rho^{univ}$ by $f_{{\mathfrak p}}$ satisfies the
three desired conditions above.  Since $\OO$ contains a field of
rationality
of $\tau$, if $\sigma \in G_K$ then ${}^{\sigma} \tau$ and $\tau$ are 
equivalent,
and so the definition of type $(k,\tau)$ is independent of the choice
of $f_{{\mathfrak p}}$.  We define
$$ R(\rhobar, k, \tau)_{\OO} = R^{univ}_{\OO} / 
\underset{{\mathfrak p} \ \text{type} \ (k,\tau)}{\bigcap} {\mathfrak p} \,. $$  
When $W(\kk)$ contains a field of
rationality
of $\tau$, we will often write 
$R(\rhobar, k, \tau)$ for $R(\rhobar,k,\tau)_{W(\kk)}$; we remark in
particular that this is always
the case for $\tau = \omt^m \oplus \omt^{lm}$,
which is rational over $\Ql$.

In the case $k=2$, if there is a surjection $\OO[[X]] \twoheadrightarrow
R(\overline{\rho},2,\tau)_{\OO}$ we say that $\tau$ is \emph{weakly
acceptable} for $\rhobar$.  If $\tau$ is weakly acceptable for $\rhobar$
and $R(\rhobar,2,\tau)_{\OO} \neq (0)$, we say that $\tau$ is acceptable
for $\rhobar$.  The above deformation rings are of
particular interest because \cite{CDT} and 
\cite{BCDT} use the methods of 
\cite{Wiles} and \cite{TaylorWiles} to prove results of roughly the
following form (for a precise statement, see Theorem 1.4.1 of 
\cite{BCDT}): 
if $\rho$ is an $l$-adic
representation of $G_{\Q}$ such that $\rho \,|_{G_{\Ql}}$ is potentially 
semistable with Galois type~$\tau$ and Hodge-Tate weights~$(0,1)$, such 
that~$\tau$ is acceptable for~$\rhobar$, and with~$\overline{\rho}$ 
modular, then $\rho$ is modular.  We remark (Appendix A of \cite{CDT}) 
that 
$$ R(\overline{\rho},k,\tau)_{\OO'} \cong \OO' \otimes_{\OO}
R(\overline{\rho},k,\tau)_{\OO} \,, $$ so when $\tau$ is defined over
the fraction field of $W(\kk)$ the
acceptability and weak acceptability of $\tau$
for $\rhobar$ depend only on $R(\rhobar,2,\tau)$.

In this article, we will prove the following cases of Conjecture 
1.2.3 of \cite{CDT}:

\begin{thm}  \label{main} Suppose that $\tau = \omt^{m} \oplus \omt^{lm}$,
where $m \in
\Z/(l^2-1)\Z$ and $m = (l+1)j+i$ with $i = 1,\ldots,l$ and $j \in
\Z/(l-1)\Z$.  Suppose also that $\rhobar \, |_{G_{\Ql}} : G_{\Ql} \rightarrow
\GL(\Fl)$, has centralizer $\Fl$, and is reducible.   Then $R(\rhobar,2,\tau)
 \neq
(0)$ only if $\rhobar
\, |_{I_l}$ has one of the following forms:
\begin{equation}
\label{eqA}
\begin{split}
\rhobar |_{I_l} & =
\begin{pmatrix}
 \omega^{i+j} & * \\
    0         & \omega^{1+j}
\end{pmatrix} \ \text{and if $i=2$, $*$ is peu-ramifi\'{e},} \\
\rhobar |_{I_l} & = 
\begin{pmatrix}
 \omega^{1+j} & * \\
    0         & \omega^{i+j}
\end{pmatrix} \ \text{and if $i=l-1$, $*$ is peu-ramifi\'{e}.}
\end{split}
\end{equation}
In each of these cases, $\tau$ is weakly acceptable for $\rhobar$.
\end{thm}

Combining \ref{main} with Theorem 1.4.1 of \cite{BCDT}, we obtain the 
following theorem:

\begin{thm}  \label{whee} Let $l$ be an odd prime, $K$ a finite extension 
of $\Ql$
in $\Qlbar$, and $\kk$ the residue field of $K$.  Let
$$ \rho : G_{\Q} \rightarrow \GL(K) $$ be an odd continuous representation
ramified at only finitely many primes.  Assume that its reduction
$$ \rhobar : G_{\Q} \rightarrow \GL(\kk) $$ is absolutely irreducible 
after
restriction to $\Q(\sqrt{(-1)^{(l-1)/2} l})$ and is modular.  Further,
suppose that
\begin{itemize}
\item $\rhobar \, |_{G_{\Ql}}$ has centralizer $\kk$,

\item $\rhobar \, |_{G_{\Ql}}$ has $\Fl$ as a field of defintion and is 
reducible,

\item $\rho \, |_{G_{\Ql}}$ is potentially Barsotti-Tate, and the
associated Weil-Deligne representation $WD(\rho \, |_{G_{\Ql}})$ is 
irreducible and tamely ramified.
\end{itemize}  
Then $\rho$ is modular.
\end{thm}

\begin{proof}  If $WD(\rho \, |_{G_{\Ql}})$ satisfies the given 
hypotheses, 
then the
$l$-type of $\rho$ is $\tau = \omt^{m}
\oplus \omt^{lm}$ for some $m$ not divisible by $l+1$.  The hypotheses on
$\rho$ guarantee that $\rhobar$ satisfies the conditions of Theorem
\ref{main}, and the very existence of $\rho$ implies that $R(\rhobar,2,\tau)
 \neq
(0)$.  Hence $\tau$ is weakly acceptable for $\rhobar$, and $\rhobar$ is
of one of the forms \eqref{eqA}.  
Once we see that (in the terminology of \cite{BCDT}) our $\tau = \omt^{m} 
\oplus \omt^{lm}$ \textit{admits} each of the two possibilities for 
$\rhobar$, then by Theorem 1.4.1 of \cite{BCDT} we obtain that $\rho$ is 
modular. 

To verify the admittance statement, one first checks that (in the notation 
of \cite{CDT} and \cite{BCDT}) $\sigma_{\tau} \cong \Theta(\chi)$ where 
$\chi : \Fllx \rightarrow \Qlbar$ maps $c \mapsto c^{-m}$.  
(The reason for an exponent of $-m$ 
instead of an exponent of $m$ is the choice of normalization
for the local Langlands correspondence in \cite{CDT}: namely, 
in Lemma 4.2.4(3) of 
\cite{CDT}, we note that since $\eta_{l,2} = \omt^{-1}$, the character
$\omt$ corresponds to $c \mapsto c^{-1}$.)  Since $m = i + (l+1)j $ with
$i \in \{1,\ldots,l\}$ and $j \in \Z/(l-1)\Z$, we may similarly write
$ -m = (l + 1 - i) + (l+1)(-1-j)$.  By Lemma 3.1.1 of \cite{CDT}, 
$\sigma_{\tau} \otimes \Flbar$ contains as Jordan-H\"older subquotients
(again, in the notation of \cite{BCDT}) the representation 
$\sigma_{l-1-i,-j}$ (if $i \neq l$) and $\sigma_{i-2,l-i-j}$ (if $i \neq 1$).
From the defintions in Section 1.3 of \cite{BCDT}, $\sigma_{l-1-i,-j}$
admits 
$\begin{pmatrix}
 \omega^{1+j} & * \\
    0         & \omega^{i+j}
\end{pmatrix}$
with $*$ peu-ramifi\'e if $i=l-1$, while $\sigma_{i-2,l-i-j}$ admits
$\begin{pmatrix}
 \omega^{i+j} & * \\
    0         & \omega^{1+j}
\end{pmatrix}$
with $*$ peu-ramifi\'e if $i=2$, as desired.
\end{proof}

\begin{remark} {\rm  Once a theory of Breuil modules with coefficients 
(see Section \ref{bm}) is sufficiently well developed, 
it should allow one to remove from Theorem \ref{whee} the 
hyptheses that $\rhobar$ is a representation defined over $\Fl$ (instead of 
over an arbitrary finite field of characteristic $l$).  One
should also then be able to use our methods to address 
Conjecture 1.2.3 of \cite{CDT} in the 
case of irreducible $\rhobar \, |_{G_{\Ql}}$.}
\end{remark}

\begin{example} {\rm Let $C$ be the genus $4$
curve $$ y^2 + (x^3 + x^2 + 1)y = -x^5 - x^4 - 2x^3 - 4x^2
- 2x - 1 \,,$$
and let $J = {\rm Jac}(C)$.  In \cite{brumer}, A. Brumer gave families of 
curves with real multiplication by $\sqrt{5}$ over $\Q$, including the family
$$ y^2 +(x^3+x+1+c(x^2+x))y=b+(1+3b)x+(1-bd+3b)x^2+(b-2bd-d)x^3-bdx^4 \,.$$
Setting $b=c=d=-1$ and substituting $y = y' + x^2$ yields the curve $C$.
Hence  $J$ carries real multiplication by 
$\sqrt{5}$, and the Galois representation 
on the $5$-adic Tate module of $J$ may be regarded as a 
two-dimensional representation $\rho_{J,5} : G_{\Q} \rightarrow 
\GL(\Q_5(\sqrt{5}))$.  In computations performed jointly with W. Stein, 
we verify that $\rho_{J,5}$ satisfies the hypotheses of Theorem 
\ref{whee}, and so $J$ is modular.   Independently, E. Gonz\'alez-Jim\'enez 
and  J. Gonz\'alez \cite{genus2} have shown the existence of a
nonconstant map $X_1(175) \rightarrow C$, and so $J$ is also modular for
that reason. }
\end{example}

The remainder of this article is concerned with the proof of Theorem
\ref{main}.

\section*{Acknowledgements}

The author is deeply indebted to his thesis supervisor, Richard Taylor.  
At every step of the way, Taylor's meticulous advice has been of 
inestimable value, and his mathematical instruction has provided a 
wonderful education.  

Brian Conrad has been a tireless resource, unfailingly willing to provide
a detailed answer to any question.  Several conversations with Fred
Diamond and Christophe Breuil were of great help at key moments.

The author is supported by an NSERC postdoctoral fellowship, and this
research was partially conducted by the author for the Clay Mathematics
Institute.

\section{Deformation theory} \label{defthy}

Henceforth $\rho$ and $\rhobar$ will denote representations of $G_{\Ql}$.
All group schemes in this article are commutative.

\subsection{Weil-Deligne representations: the Barsotti-Tate case} 
\label{btcase}
When $\rho : G_{\Ql} \rightarrow {\rm GL}_d(K)$ is potentially Barsotti-Tate, 
we provide an alternate description of $WD(\rho)$, directly following 
Appendix B.3 of \cite{CDT}.
Suppose $\rho$
becomes Barsotti-Tate over a finite Galois extension $E$ of $\Ql$, 
so that $\rho \,|_{G_E}$ arises from
an $l$-divisible group $\Gamma$ over $\OO_E$.  Write $\OO$ for the integers
of $K$, and  $\kk$ for the residue field of $E$.

By Tate's full faithfulness theorem (Theorem 4 of \cite{tate}), 
$\Gamma$ has an action of $\Gal(E/\Ql)$ over the action of $\Gal(E/\Ql)$
on $\Spec(\OO_E)$.  This reduces to an action on the closed fibre $\Gamma \times \kk$.  Let $\phi$ be the Frobenius endomorphism of the closed fibre
of $\Gamma$; then we produce an action of $W_l$ on $\Gamma_{/\kk}$ 
by letting $g$ act via $g \,|_E \circ \phi^{-v(g)}$.  

This above action of $W_l$ is a right-action.  It therefore translates into
a left-action on the contravariant Dieudonn\'{e} module $D(\Gamma_{/\kk})$.
Then $D$ is a free $W(\kk)$-module of rank $d[K : \Ql]$.
Let $F$  denote the Frobenius element of the Dieudonn\'{e} ring.

Next, we define an action of $W_l$ on
\begin{equation}\label{eqH}
D'(\Gamma_{/\kk}) = \Hom_{W(\kk)} (D(\Gamma_{/\kk}), W(\kk)) \,.
\end{equation}
We set
$$ \phi'(f) = \sigma \circ f \circ F^{-1} $$
on $ D'(\Gamma_{/\kk})[1/l]$, where $\sigma$ is Frobenius on $W(\kk)$, and for $g \in \Gal(E/\Ql)$ we set
$$ g(f) = \overline{g} \circ f \circ g^{-1} $$ 
where $\overline{g}$ is the map $g$ induces on $W(\kk)$ and $g^{-1}$ 
is the semilinear action on $D(\Gamma_{/\kk})$ coming from the semilinear
action on $\Gamma$.  Finally, as usual, we let $W_l$ act on $D'(\Gamma_{/\kk})$ by letting $g$ act as $g \,|_{E} \circ (\phi')^{-v(g)}$.

Finally, we note that the action of $\OO$ on $\Gamma$ propagates through
all of the above constructions, and we have (Proposition B.3.1 in 
\cite{CDT}): 
$$ WD(\rho) \cong D'(\Gamma_{/\kk}) \otimes_{W(\kk) \otimes_{\Zl} \OO} \Qlbar \,.$$

\subsection{Dieudonn\'{e} module calculations} \label{dieudonne} 

For the rest of this article, we fix $\tau =
\omt^m \oplus \omt^{lm}$, 
and the following notation.  Let~$\Qll$ be the copy in $\Qlbar$ of
the field of fractions of the Witt vectors $W(\Fll)$, and let $\pi$ be a
choice of $\unif$.  Let $E = \Ql(\pi)$, $E' = \Qll(\pi)$.  Note that $\tau 
\,|_{I_E}$ is trivial.  We will regard an element $\zeta \in \Fll$ as an 
element in $W(\Fll)$  (and hence in $\Qll$) via the Teichmuller lifting 
map.  Let $g_{\zeta}$ denote the element of $\Gal(E'/\Ql)$ fixing~$\Qll$ 
and sending $\pi$ to $\zeta \pi$.  Let $\varphi$ denote the element of 
$\Gal(E'/\Ql)$ fixing $\pi$ and extending the nontrivial automorphism of 
$\Qll$.  

Suppose $\rho : G_{\Ql}
\rightarrow \GL(K)$ is a potentially Barsotti-Tate representation with
$l$-type $\tau$ and with determinant $$\det(\rho)=
\epsilon \cdot {\rm Teich}(\omega^{-1} \det(\rhobar)) \,,$$ where Teich
denotes the Teichm\"{u}ller lift.   

We now specialize the discussion of Section \ref{btcase} to 
this
 situation.
We know $\tau \, |_{I_{E'}}$ is trivial, and so $\rho$
becomes Barsotti-Tate when restricted to $G_{E'}$.  Consequently, we
obtain an $l$-divisible group $\Gamma$ over $\OO_{E'}$ such that the Tate
module of the generic fibre of $\Gamma$ is $\rho \, |_{G_{E'}}$.  
The field residue field of $E'$ is $\kk = \Fll$, the Witt vectors $W(\kk)=\Zll$, $\sigma$ is the Frobenius automorphism of $\Zll$, the map $\overline{g}_{\zeta}$ is the identity for each $\zeta$, and the map $\overline{\varphi}=\sigma$.

We saw in Section \ref{btcase} (Proposition B.3.1 of \cite{CDT})
that 
\begin{equation}\label{eqB}
WD(\rho) \cong D'(\Gamma_{/\Fll}) \otimes_{\Zll \otimes_{\Zl} \OO_K} \Qlbar \,,
\end{equation}
where $g \in W_l$ acts on the right-hand side via 
$g \, |_{E'} \circ (\phi')^{-v(g)}$.
In particular, $I_l$
acts via $I_l \twoheadrightarrow \Gal(E'/\Qll)$, and since $v(I_l)=0$,
no untwisting is needed.

Since $\tau = \omt^{m} \oplus \omt^{lm}$, there exist
basis elements $\mathbf{v},\mathbf{w}$ of $D'(\Gamma_{/\Fll}) \otimes_{\Zll
\otimes_{\Zl}
\OO_K} \Qlbar$ so that for $g \in I_l$,
$$ g(\mathbf{v}) = \omt^{m}(g) \mathbf{v} $$
and $$ g(\mathbf{w}) = \omt^{lm}(g) \mathbf{w}. $$

For $\zeta \in \Fll$, by definition we have
$$ \omt^m(g_{\zeta}) = (g_{\zeta}(\pi)/\pi)^m = \zeta^m \,,$$
and similarly
$\omt^{lm}(g_{\zeta}) = \zeta^{lm}$.  
Thus $g_{\zeta}(\mathbf{v}) = \zeta^{m} \mathbf{v}$
and $g_{\zeta}(\mathbf{w}) = \zeta^{lm} \mathbf{w}$ .  Similarly, we find
$g_{\zeta}^l(\mathbf{v}) = \zeta^{lm} \mathbf{v}$ and 
$g_{\zeta}^l(\mathbf{w}) = \zeta^{m}
\mathbf{w}$, from which we conclude that $g_{\zeta} + g_{\zeta}^l$ acts on 
$D'(\Gamma_{/\Fll}) \otimes_{\Zll \otimes_{\Zl} \OO_K} \Qlbar$ by scalar
multiplication by $\zeta^{m} + \zeta^{lm}$, whereas $g_{\zeta}^{l+1}$ acts by
scalar multiplication by $\zeta^{(l+1)m}$.  (The action is linear, and not
semilinear, since the image of $I_l \rightarrow \Gal(E'/\Qll)$ acts
trivially on the coefficients $\Zll$.)

We now wish to use \eqref{eqH} and the 
action $g_{\zeta}(f) = \overline{g}_{\zeta} 
\circ f \circ g_{\zeta}^{-1} = f \circ g_{\zeta}^{-1}$ on 
$D'(\Gamma_{/\Fll})$ to understand the action of $g_{\zeta}$
on $D(\Gamma_{/\Fll})$.

Since $D'(\Gamma_{/\Fll})$ is a free module, the actions of $g_{\zeta}$ and
$g_{\zeta}^l$ must sum and multiply on $D'(\Gamma_{/\Fll})$ to scalar 
multiplication
by $\zeta^{m}+\zeta^{lm}$ and $\zeta^{(l+1)m}$ respectively.  If $f \in
D'(\Gamma_{/\Fll})$, we know 
$$g_{\zeta}(f(x)) = f(g_{\zeta}^{-1} x)$$ with $x \in D(\Gamma_{/\Fll})$.
It follows that $$ f((g_{\zeta}^{-1} + g_{\zeta}^{-l})x) =
(g_{\zeta} + g_{\zeta}^{l}) f(x) = (\zeta^{m} + \zeta^{lm}) f(x)
= f((\zeta^{m} + \zeta^{lm}) x). $$  
By freeness, for any nonzero $x \in D(\Gamma_{/\Fll})$
we can find $f \in D'(\Gamma_{/\Fll})$ which does not vanish on $x$, so we
conclude
that $g_{\zeta}^{-1} + g_{\zeta}^{-l}$ acts as scalar multiplication by
$\zeta^{m} + \zeta^{lm}$ on $D(\Gamma_{/\Fll})$.  Replacing $\zeta$
by $\zeta^{-1}$, we have found that 
\begin{equation}\label{eqC}
g_{\zeta} + g_{\zeta}^l \ \text{acts as scalar
multiplication by} \ \zeta^{-m} + \zeta^{-lm} \ \text{on} \ 
D(\Gamma_{/\Fll}) \,.
\end{equation}
Similarly 
\begin{equation}\label{eqD}
g_{\zeta}^{l+1} \ \text{acts as scalar multiplication by} \ 
\zeta^{-(l+1)m} \ \text{on} \ D(\Gamma_{/\Fll}) \,.
\end{equation}
We next wish to see what the determinant condition tells us about
$D(\Gamma_{/\Fll})$. 
Let $\chi_l$ denotes the
$1$-dimensional unramified character of $W_l$ sending arithmetic 
Frobenius to $l$, and let $$\chi =
{\rm Teich}(\omega^{-1} \det(\rhobar)) \, |_{W_{\Ql}} \otimes_{K}
\Qlbar \,.$$
By the examples in Section B.2 of \cite{CDT}, and since $WD$ is compatible 
with tensor products,
$$WD(\det(\rho)) = \chi_l \chi \,.$$  Let $s$ be any lift of $\varphi$
to $W_{l}$, so $s$ is a lift of arithmetic Frobenius but fixes $F$.  
Since $WD$ is
compatible with the
formation of exterior products, we know $\det(WD(\rho)) = WD(\det(\rho))$, 
and in particular $\det(WD(\rho)(s)) = l T$, where $T =
{\rm Teich}(\omega^{-1} \det(\rhobar))(s) = {\rm
Teich}(\det(\rhobar))(s)$.  
(We have $\omega(s)=1$ since $s$ fixes $F$.)
Note that $T$ depends only on $\rhobar$, not on $\rho$ or the choice of
$s$.

We claim that $\Trace(WD(\rho)(s)) = 0$.  Since
$$ WD(\rho) \,|_{I_l} =
\begin{pmatrix}
\omt^{m} & 0 \\
0 & \omt^{lm}
\end{pmatrix}
$$
and since for any $u \in I_l$ we have the relation $WD(\rho)(sus^{-1})
= WD(\rho)(u^l)$, it follows immediately that $WD(\rho)(s)$ 
must act via a matrix
$$\begin{pmatrix}   
    0  &  * \\
    *  &  0
\end{pmatrix}.$$
Therefore, we have shown that $WD(\rho)(s)$ satisfies the
characteristic polynomial $X^2 + lT = 0$.
By \eqref{eqB}, and since 
$D'(\Gamma_{/\Fll})$ is free, the action of $s$ on 
$D'(\Gamma_{/\Fll})$
must
satisfy this same 
polynomial.

For $D(\Gamma_{/\Fll})$, note that 
if $f \in D'(\Gamma_{/\Fll})[1/l]$ then 
$$s(f)(x) = \varphi \circ (\phi')^{-1} (f)(x) = (\sigma (\sigma^{-1}
\circ f \circ F) \circ \varphi)(x) = f( F \circ \varphi (x))$$ for $x \in
D(\Gamma)$.  Then $$ s^2 (f)(x) = f (F^2 \circ \fr^2
(x)) = f(F^2 (x)).$$  Since $s^2 + lT = 0$ on $D'(\Gamma_{/\Fll})$,
we learn that $$ f( (F^2 + lT)x ) = 0 $$ for all $x$ and $f$,
and
consequently $F^2 + lT = 0$ on $D(\Gamma_{/\Fll})$.  Applying $V$ to both sides of
this equation we see $F(FV) + lTV = lF + lTV = 0$, and since $D(\Gamma_{/\Fll})$
is free we obtain the relation 
\begin{equation}\label{eqE}
F + TV = 0
\end{equation}
on $D(\Gamma_{/\Fll})$.
 
\subsection{Deformation problems}  We will make use of the following 
definitions, essentially following Section 4 of \cite{BCDT}:

\begin{defn} 
{\rm
If $X$ is a scheme over $\Spec A$ and $g: A \rightarrow B$ is a 
ring homomorphism, let $^gX$ denote the pullback of $X$ by $g$.
Suppose that $K/L$ is an Galois extension of fields over $\Ql$, and let
$\G$ be a group scheme over $\OO_{K}$.  By \textit{generic 
fibre descent data} from $K$ to $L$, we mean a collection of isomorphisms
$$ [g] : \G \rightarrow \ ^g \G $$
for $g \in \Gal(K/L)$
satisfying the compatibility conditions $[gh]=(^g[h])\circ[g]$. 
The pair $(\G, \{[g]\})$, which we will sometimes abbreviate as $\G$,
is a \textit{group scheme with descent data}.
Note that since $\OO_K/\OO_L$ is not necessarily \'etale, we do not 
necessarily obtain a descended group scheme over $\OO_L$.  However,
since $K/L$ is \'etale we can descent the generic fibre as usual, 
and we denote the descended
$L$-group scheme by $(\G,\{[g]\})_{L}$.  By the descended $G_L$-representation
of $(\G,\{[g]\})_{L}$, we mean the representation of $G_L$ on
$(\G,\{[g]\})_{L}(\Qlbar)$.

}
\end{defn} 

\begin{defn} {\rm
If $G$ is a finite
group scheme over a field $L/\Ql$, then an \textit{integral model}
of $G$ is a finite flat group 
scheme $\G$ over $\OO_L$ such that $\G \times_{\OO_L} L
\cong G$.  More generally, if $K/L$ is a Galois extension, then an 
\textit{integral model of $G$ with descent data} over $\OO_K$ is a  
finite flat group
scheme $(\G,\{[g]\})$ over $\OO_K$ with descent data to $L$ such 
that $(\G,\{g]\})_L \cong G$.
}
\end{defn}

Fix $\rhobar : G_{\Ql} \rightarrow \GL(\Fl)$ a reducible 
Galois representation
with
centralizer $\Fl$, let $M_{\rhobar}$ denote the standard 
$\Fl$-vector space on which $G_{\Ql}$
acts via $\rhobar$, and let $T = {\rm Teich}(\det(\rhobar))(\Fr)$. 

We let
$\s(\rhobar)$ denote the full subcategory of the category of finite length
discrete $\Zl$-modules with $\Zl$-linear action of $G_{\Ql}$ consisting
of objects which admit a finite filtration such that each graded piece 
is isomorphic to $M_{\rhobar}$.  
Let $\s$ be the full subcategory of $\s(\rhobar)$
consisting of objects $X$ for which there exists a finite flat
$\OO_{E'}$-group scheme $(\G,\{[g]\})$ with descent data to $\Ql$
 such that $X \cong
(\G,\{[g]\})_{\Ql}(\Qlbar)$ as $\Zl[G_{\Ql}]$-modules and such 
that $\z + \z^l
- (\zeta^{-m} + \zeta^{-lm})$ and $\z^{l+1} - \zeta^{-(l+1)m}$ for all
$\zeta \in \Fllx$, as well as $F
+ TV$, annihilate the Dieudonn\'{e} module $D(\G \times \Fll)$.

From  Lemma 4.1.3 of \cite{BCDT} it follows that $\s$ is closed under finite
products, subobjects, and quotients.   Following Section 4.3 of 
\cite{BCDT}, define the set-valued functor
${\mathcal D}^{\s}_{\rhobar,\Zl}$ on the category of complete Noetherian local 
$\Zl$-algebras with residue field $\Fl$ by letting ${\mathcal
D}^{\s}_{\rhobar,\Zl}(R)$ be the set of conjugacy classes of continuous
$R$-representations such that $\rho \mod \mathfrak{m}_{R}$ is conjugate to
$\rhobar$ and such that for each open ideal $\mathfrak{a} \subset R$ the
action of $\rho$ makes $(R/\mathfrak{a})^2$ into an object of $\s$.

By a theorem of Ramakrishna \cite{ram}, if ${\mathcal
D}^{\s}_{\rhobar,\Zl}(\Fl)$ is nonempty, then the functor ${\mathcal
D}^{\s}_{\rhobar,\Zl}$ is representable; in this case, let 
$R^{\s}_{\rhobar,\Zl}$
denote the resulting deformation ring.  We have:

\begin{prop}  If  ${\mathcal D}^{\s}_{\rhobar,\Zl}(\Fl)$ is 
nonempty, then there
is a surjection $$R^{\s}_{\rhobar,\Zl} \twoheadrightarrow \Rdvz \,.$$
\end{prop}

\begin{proof}  Let $R^{univ}_{\rhobar}$ denote the universal deformation ring
for $\rhobar$, and let $${\mathcal I} = \ker(R^{univ}_{\rhobar}
\twoheadrightarrow R^{\s}_{\rhobar,\Zl}) \,.$$  It suffices 
to show $${\mathcal I}
\subset \ker(R^{univ}_{\rhobar} \twoheadrightarrow \Rdvz) =
\underset{\mathfrak{p}
\, \text{type} \, \tau}{\bigcap} \mathfrak{p} =  \underset{i \ge 1,
\mathfrak{p} \, \text{type} \, \tau}{\bigcap} (\mathfrak{p},l^i).$$  
In other words, we need to show that each map $R^{univ}_{\rhobar}
\twoheadrightarrow R^{univ}_{\rhobar} / (\mathfrak{p},l^i)$ factors through
$R^{univ}_{\rhobar} \twoheadrightarrow R^{\s}_{\rhobar,\Zl}$.  
Let $\tilde{\rho}$ 
denote the representation arising from  $R^{univ}_{\rhobar} 
\twoheadrightarrow
R^{univ}_{\rhobar} / (\mathfrak{p},l^i)$.  
Since $\mathfrak{p}$ has type $\tau$,
there is an extension $K/\Ql$ and an exact sequence $0 \rightarrow
\mathfrak{p} \rightarrow R^{univ}_{\rhobar} \rightarrow K$ so that the 
resulting
$\rho : G_{\Ql} \rightarrow \GL(K)$ is of type~$\tau$.  The results of
Section \ref{dieudonne} produce an $l$-divisible group $\Gamma/\OO_{E'}$
satisfying the desired relations on the Dieudonn\'{e} module of its closed
fibre and whose generic fibre representation  is $\rho \, |_{G_{E'}}$.
The $l^i$-torsion $\Gamma[l^i]$ is the desired finite flat group scheme with
descent data which shows that the conjugacy class of $\tilde{\rho}$ is
indeed in  ${\mathcal D}^{\s}_{\rhobar,\Zl}(R^{univ}_{\rhobar} 
/ (\mathfrak{p},l^i))$.
\end{proof}

\begin{cor} If  ${\mathcal D}^{\s}_{\rhobar,\Zl}(\Fl)$ is nonempty, then the
dimension of the tangent space of $\Rdvz$ is at most
the dimension of the tangent space of $R^{\s}_{\rhobar,\Zl}$.
\end{cor}

The rest of this article will be concerned with the proof of 
the following theorem:

\begin{thm} \label{rest} If $\Rdvz \neq (0)$, then ${\mathcal
D}^{\s}_{\rhobar,\Zl}(\Fl)$
is nonempty.  In this case, $\rhobar \, |_{I_l}$ is of one of the two
forms \eqref{eqA}, 
and up to isomorphism there is exactly one finite
flat group scheme $(\G,\{[g]\})$ over $\OO_{E'}$ with descent data to
$\Ql$
such that $(\G,\{[g]\})_{\Ql} \cong \rhobar$ and such that $D(\G \times
\Fll)$ satisfies the relations (\ref{eqC}, \ref{eqD}, \ref{eqE}).
The space of extensions
of  $(\G,\{[g]\})$ by  $(\G,\{[g]\})$, in the category of finite flat
$\OO_{E'}$-group schemes
with descent data,
whose Dieudonn\'{e} modules satisfy these relations, is
$1$-dimensional.
\end{thm}

This theorem evidently implies that if $\Rdvz \neq (0)$ then
$R^{\s}_{\rhobar,\Zl}$ exists and has a $1$-dimensional tangent space, which 
completes the proof of Theorem \ref{main}.

\subsection{Strategy of the calculation}

If $\Rdvz \neq (0)$, then there exists a prime $\mathfrak{p}$ of type
$\tau$.  Hence there is a lift $\rho$ of $\rhobar$ which arises from an
$l$-divisible group $\Gamma$ over $\OO_{E'}$ with descent data to $\Ql$ and
satisfying the Dieudonn\'{e} module relations (\ref{eqC}, \ref{eqD}, 
\ref{eqE}).  
Then the
$l$-torsion $\Gamma[l]$ is filtered by integral models 
for $\rhobar$ with descent data, 
so we see that  ${\mathcal D}^{\s}_{\rhobar,\Zl}(\Fl)$
is nonempty.

It remains: to determine all (reducible) $\rhobar$ for which there exists
a group scheme $(\G,\{[g]\})$ over $\OO_{E'}$ with descent data to $\Ql$,
such that $(\G,\{[g]\})_{\Ql} \cong \rhobar$, and such that the 
Dieudonn\'{e} module of $\G$ satisfies the relations (\ref{eqC}, \ref{eqD}, 
\ref{eqE}); to show that when such a group scheme with descent data
exists, there is exactly one of them; and, in this case, to compute the
extensions described in Theorem \ref{rest}.

Note that since $\rhobar$ is reducible, any integral model
with descent data for $\rhobar$ over $\OO_{E'}$
is an extension of rank $1$ group schemes which are integral models
for the sub- and quotient- characters of $\rhobar \, |_{G_{E'}}$,
and by a 
 scheme-theoretic
closure argument (see \cite{Raynaud} or Lemma 4.1.3 of \cite{BCDT})
this is actually an extension in the category of integral models with
descent data.  In Section \ref{bm}, we describe the category 
of Breuil modules with descent data, 
which is anti-equivalent to the category of 
finite flat $l$-torsion group schemes with descent data.
We will proceed to use Breuil's theory
as follows.  We
compute all of the rank $1$ Breuil modules with descent data
over a tamely ramified extension, and identify explicitly the
Galois characters to which these Breuil modules with descent data 
correspond.
We next classify all of the (rank $2$) extensions, in the
category of Breuil modules with descent data.
of these rank $1$ Breuil modules with descent data, after 
which we may discard from 
consideration those Breuil modules with descent data
 which do not correspond to group
schemes satisfying the relations (\ref{eqC}, \ref{eqD}, \ref{eqE}) on 
their Dieudonn\'{e} modules.

We will indeed see that the only $\rhobar$ for which
integral models exist that admit generic fibre descent data satisfying the 
desired 
Dieudonn\'{e} module relations, are exactly those of the form \eqref{eqA}.
Moreover, for each such $\rhobar$
this integral model with descent data 
will be seen to be unique (up to isomorphism).  
We complete the proof of Theorem \ref{rest} by
calculating $\textrm{Ext}^1(\M,\M)$, 
in the category of Breuil modules with 
descent data, for the Breuil modules with descent data $\M$ 
corresponding to these integral models with descent data,
and checking that the space of extensions satisfying the Dieudonn\'{e} 
module relations is at most $1$-dimensional.

\section{Review of Breuil modules with descent data} \label{bm}

\subsection{Breuil modules}

We remind the reader that the prime $l$ is odd. 
Let $K/\Ql$ be a finite extension, and suppose $K$ has integers $\OO$,
ramification index $e_K$, and residue field $\kk$.  Fix a uniformizer
$\pi$ of $\OO$.  A Breuil module $(\M,\MM,\p)$ for $K$
consists of the following data:
\begin{itemize}
\item a finite-rank free $\kk[u]/u^{e_K l}$-module $\M$,
\item a submodule $\MM \subset \M$ such that $\MM \supset u^{e_K} \M$, and
\item an additive map $\p: \MM \rightarrow \M$ such that $\p(hv)=h^l \p(v)$
for any $h \in \kk[u]/u^{e_Kl}$ and $v \in \MM$, and such that the 
$\kk[u]/u^{e_K l}$-span of $\p(\MM)$ is all of $\M$.
\end{itemize}
Morphisms of Breuil modules are $\kk[u]/u^{e_K l}$-module
homomorphisms which
preserve $\MM$ and commute with $\p$.  The rank
of a Breuil module is defined to be its rank as a $\kk[u]/u^{e_K l}$-module.

C. Breuil \cite{Br,Breuil0} has proved the following theorem:

\begin{thm}  There is an (additive) contravariant
equivalence of categories,
depending on the choice of uniformizer $\pi$, between the
category of Breuil modules for $K$ and the category of finite flat
group schemes over $\OO$ which are killed by~$l$.  The 
rank of the Breuil module is the same as the rank of the corresponding
group scheme.
\end{thm}

When the field $K$ and uniformizer $\pi$ are clear from context,
by the Breuil module corresponding to a group scheme we will mean 
the Breuil module obtained from the group scheme 
via this equivalence (and vice-versa).

The Breuil module functor has numerous useful properties: for
example, a short exact sequence of group schemes maps under the functor to
a short exact sequence of Breuil modules, and a sequence
of Breuil modules is short-exact if and only if the
underlying sequence of 
$\kk[u]/u^{e_K l}$-modules is short-exact. 
(\cite{BCDT},
Lemma 5.1.1.)  This will allow
us directly to compute Exts of Breuil modules. 

\vskip 0.2cm

There is a very useful compatibility between Breuil theory and
contravariant Dieudonn\'{e} theory.  Let $$u^{e_K} - l G_{\pi}(u)$$
be the minimal polynomial of $\pi$ over $W(\kk)$, and let
$c_{\pi} = -G_{\pi}(u)^l \in \kk[u]/u^{e_K l}$.  
On any Breuil module, define 
$\phi : \M \rightarrow \M$ via
$$\phi(v) = \frac{1}{c_{\pi}} \p (u^{e_K} v).$$  Note that $u^{e_K} v 
\in \MM$ 
by definition.  Then (\cite{BCDT}, Theorem 5.1.3(3)) if $\M_{\pi}$ is
the Breuil module corresponding to the group scheme $\G$ (with $\pi$
as our fixed uniformizer), there is a 
canonical $\kk$-linear isomorphism
\begin{equation}\label{eqG}
D(\G) \otimes_{\kk,{\rm Frob_l}} \kk \cong \M_{\pi}/u\M_{\pi}
\end{equation}
under which $F \otimes \Frob_l$ corresponds to $\phi$ and $V \otimes
\Frob_l^{-1}$ corresponds to the composition
$$ \M/u\M \overset{\p^{-1}}{\longrightarrow} 
\MM/u\MM \rightarrow \M/u\M \,.$$
(One can see that $\p$ mod ${u}$ is always bijective.)

\subsection{Rank 1 Breuil modules}

It is an
informative exercise (\cite{BCDT}, Example 5.2)  
to check that
the rank $1$ Breuil modules are of the form:
$$ \M = (\kk[u]/u^{e_K l}) \e \,,$$ 
$$ \MM = (\kk[u]/u^{e_K l}) u^r \e \,,$$
$$ \p(u^r \e) = a \e $$
with $0 \le r \le e$ and $a \in \kk^{\times}$.  We will denote this
module as $\M(r,a)$.   We recommend that
the reader verify that nonzero homomorphisms $\M(r,a) \rightarrow \M(r_1,a_1)$
exist if and only if $r_1 \ge r$, $r_1 \equiv r \pmod{l-1}$, and
$a/a_1 \in (\kk^{\times})^{l-1}$, and are
given exactly by the linear maps
$\e \mapsto bu^{l(r_1-r)/(l-1)} \e_1$ where $b^{l-1}=a/a_1$.

From 3.1.2 of \cite{Br}, the affine algebra underlying the group scheme
corresponding to $\M(r,a)$ is $$\OO[X]/(X^l + \frac{\pi^{e_K -r}\tilde{a}}
{G_{\pi}(\pi)} X) \,,$$ where $\tilde{a}$ denotes the Teichm\"uller lift of
$a$.  We note that we may say even more, namely that the comultiplication
on this algebra is exactly that which one would expect from the Oort-Tate
classification, namely:

\begin{lemma} \label{comult} 
Set $C = \frac{\pi^{e_K-r} \tilde{a}}{G_{\pi}(\pi)}$.
The group scheme corresponding to $\M(r,a)$ (under
the fixed choice of uniformizer $\pi$) is isomorphic to the group scheme 
with affine algebra 
$\OO[X]/(X^l + C X)$
and comultiplication
\begin{equation}\label{eqF}
 X \mapsto 1 \otimes X + X \otimes 1 - \frac{l}{C}
\sum_{i=1}^{l-1} \frac{X^i}{w_i} \otimes \frac{X^{l-i}}{w_{l-i}}
\end{equation}
where the units $w_i \in \Zlx$ are as defined in Section 2 of 
\cite{oorttate}.
\end{lemma}

\begin{proof}  If $r < e_K$, so that $C$ is divisible by $\pi$,
one simply needs to note that
$$ \OO[X]/(X^l + CX) \cong \OO[X]/(X^l + C'X) $$ if and only if $C/C' \in 
(\OO^{\times})^{l-1}$.  By the Oort-Tate classification, the group scheme
corresponding to $\M(r,a)$ is isomorphic to some $\Spec \OO[X]/(X^l + C'X)$
with comultiplication as in \eqref{eqF} with $C'$ in place of $C$.  Since $C'/C  \in 
(\OO^{\times})^{l-1}$, it is therefore also isomorphic to  
$\Spec \OO[X]/(X^l + CX)$ with comultiplication \eqref{eqF}.

If $r=e_K$, the Dieudonn\'e module compatibility \eqref{eqG} shows that
the classical Dieudonn\'e module of the closed fibre of the group 
scheme corresponding to 
$\M(e_K,a)$ is isomorphic to 
$$ \kk[F,V]/ (F + \frac{a}{G_{\pi}(0)}, V ) \,,$$
where the ring $\kk[F,V]$ is noncommutative if $\kk \neq \Fl$, satisfying
$Fx = x^l F$, $Vx^l = xV$ for $x \in \kk$, and $FV=VF=0$.
We recall from Section 3 of \cite{oorttate} that the 
Dieudonn\'e modules attached to the group scheme
$$ \kk[X]/(X^l - \alpha X) \,, \ X \mapsto 1 \otimes X + X \otimes 1 + 
\beta \sum_{i=0}^{l-1} \frac{X^i}{i!} \otimes \frac{X^{l-i}}{(l-i)!} $$
with $\alpha \beta =0$ is $$ \kk[F,V]/ (F - \alpha, V - \beta^{1/l}) \,.$$
If the group scheme corresponding to $\M(e_K,a)$ is isomorphic to 
$\Spec \OO[X]/(X^l + CX)$ with comultiplication \eqref{eqF}, it follows that 
$C \in \OO^{\times}$, and that the image of $C$ in $\kk$ differs from 
that of $ -a/G_{\pi}(0)$ by an $(l-1)^{\textrm{st}}$ power.  Noting that 
$G_{\pi}(0)$ and $G_{\pi}(\pi)$ have the same image in $\kk$, the claim
follows using Hensel's lemma.
\end{proof}

\subsection{Generic fibre descent data}

A group scheme with descent data $(\G,\{ [g] \})$ corresponds
to a Breuil module with descent data.  In the case where $K/L$ a is
tamely ramified Galois extension with relative ramification index
$e = e(K/L)$, and the uniformizer $\pi$ of $K$
satisfies $\pi^{e} \in L$, the description of generic fibre
descent data is 
fairly
simple.  (In the wild case, it is decidedly not.)  The following
description has not appeared in the literature, but is
essentially a transcription from an unpublished preprint of 
B. Conrad \cite{WRDR}, included here by permission:

\begin{thm} If $K/L$ is a tamely ramified Galois extension 
and $\pi^{e} \in L$, then
giving generic fibre descent data on $(\M,\MM,\p)$ is equivalent to 
giving, for each
$g \in \Gal(K/L)$, an additive bijection $[g]: \M \mapsto \M$ satisfying:
\begin{itemize}
\item each $[g]$ preserves $\MM$ and commutes with $\p$,
\item $[1]$ is the identity and $[g][h] = [gh]$, and
\item $g(a u^i v) = g(a)(zu)^i g(v)$, where $g(\pi)=z\pi$ and $a \in \kk$
is regarded as being in $K$ via the Teichm\"{u}ller lift.
\end{itemize}
Moreover, this generic fibre 
descent data is compatible with Dieudonn\'e 
theory \eqref{eqG}.
\end{thm}

To see that this description follows from the (significantly) more 
involved description, found in Section 5.6
of \cite{BCDT}, of
generic fibre descent data over general (i.e., possibly wild)
field extensions, we again
quote from \cite{WRDR}: observe,
in the notation of \cite{BCDT}, that we can choose $H_g(u)=g(\pi)/\pi$
a root of unity, so $t_g = 0$ for all $g$, $f_{g_1,g_2} = 0$ for all
$g_1,g_2$, and therefore $\textbf{1}_{g_1,g_2}$ is the identity.

Given two Breuil modules $\M'$, $\M''$ with descent data, an extension in the 
category of Breuil modules with descent data is an extension of Breuil modules
$$ 0 \rightarrow \M' \rightarrow \M \rightarrow \M'' \rightarrow 0$$
with generic fibre
 descent data on $\M$ such that for all $[g]$ the following diagram 
commutes:
$$ \xymatrix{
 0 \ar[r] & \M' \ar[r] \ar[d]^{[g]} & \M \ar[r] \ar[d]^{[g]} & \M'' \ar[d]^{[g]}\ar[r] & 0 \\
 0 \ar[r] & \M' \ar[r] & \M \ar[r] & \M'' \ar[r] & 0     } $$

If $\M$ is a Breuil module with descent data corresponding to a group scheme
$\G$ with descent data from $K$ to $L$, then by the descended $L$-group
scheme (resp. $G_L$-representation) of $\M$, we mean the descended $L$-group
scheme (resp. $G_L$-representation) of $\G$.

For any further facts about Breuil modules which may be necessary,
the reader can refer to Section 5 of \cite{BCDT}.  We now begin the
computations needed for the proof of Theorem \ref{rest}.

\section{A Galois cohomology lemma}

Let $K/L$ be a tamely ramified Galois extension of $\Ql$, 
let $\kk/{\mathbf l}$ be the extension of residue 
fields, let 
$e = e(K/L)$ be the ramification index, and let $\pi \in K$ be a 
uniformizer such that $\pi^e \in L$.  For ease of notation, we will 
identify the elements of $\kk$ with their Teichm\"uller lifts in $K$.  

\begin{lemma} \label{galcoh} With the above notation, let $n$ be a 
positive integer 
and let $G = \Gal(K/L)$ act on 
$\kk[u]/u^n$ via ${}^g ( \sum_{i=0}^{n-1} a_i u^i ) = \sum_{i=0}^{n-1} 
g(a_i) 
\left(\frac{g\pi}{\pi}\right)^i u^i$.  Under this action,
\begin{itemize}
\item $\HH^q(G,\kku) = 0$ for all $q > 0$, and
\item $\#\HH^1(G,\kkux) = e$.
\end{itemize}
The nonzero elements of $\HH^1(G,\kkux)$ are represented by the cocycles 
$g \mapsto \left(\frac{g\pi}{\pi}\right)^i$, for $0 \le i < e$.
\end{lemma}

\begin{proof}
Let $\kk_{i}$ denote the additive group $\kk$ with the $G$-action $g 
\cdot a 
= 
g(a) \left(\frac{g\pi}{\pi}\right)^i$.  Let $I$ denote the inertia 
subgroup of $G$, and let $K_0 = K^I$ be the maximal unramified extension 
of $L$ inside $K$.  Then $\HH^q(I,\kk_i)=0$ for all $q > 0$ as $\#\kk_i$ 
is a 
power of $l$ whereas $\#I$ is prime-to-$l$ (since $L/K$ is tame).  As a 
result, the Hochschild-Serre spectral sequence
provides an isomorphism 
$$\HH^q(G/I,\kk_i^{I}) \xrightarrow{\sim} \HH^q(G,\kk_i)$$
for all $q$.  We compute: $a \in \kk_i^I$ if and only if $g(a) 
\left(\frac{g\pi}{\pi}\right)^i = a$ for all $g \in I$, if and only if $a 
\pi^i \in K_0$.  Noting that $a \in K_0$, if $a \neq 0$ then $\pi^i \in 
K_0$, and in this case we see by considering the valuation that that $i$ 
is divisible by $e$.  Hence either $\kk_i^{I}=0$, or $\kk_i^{I} = \kk$ 
with the usual action of $G/I=\Gal(\kk/{\mathbf l})$.  In both cases 
$\HH^q(G/I,\kk_i^{I})=0$, and so $\HH^q(G,\kk_i)=0$ for all $i$ and all 
$q > 0$.

Now the first claim of the lemma follows immediately from the isomorphism 
$\kku = \oplus_{i=0}^{n-1} \kk_i$.  To see the second claim, observe that 
we have a short exact sequence of $G$-modules
\begin{equation}\label{eqI}
 0 \rightarrow \{ 1 + a u^n + u^{n+1} \kku \} \rightarrow 
(\kk[u]/u^{n+1})^{\times} \rightarrow \kkux \rightarrow 0
\end{equation}
and note that the first module in the sequence \eqref{eqI} 
is isomorphic to $\kk_n$ 
when $n \ge 1$.  
From the long exact sequence of cohomology associated to \eqref{eqI} 
and the vanishing of 
$\HH^1(G,\kk_n)$ and $\HH^2(G,\kk_n)$ we obtain an isomorphism 
$\HH^1(G,(\kk[u]/u^{n+1})^{\times}) \cong \HH^1(G,\kkux)$ for all $n \ge 
1$.  Hence $\HH^1(G,\kku) \cong \HH^1(G,\kk^{\times})$.  By another 
application of the Hochschild-Serre spectral sequence, we obtain
\begin{equation}\label{eqJ}
0 \rightarrow \HH^1(G/I,(\kk^{\times})^I) \rightarrow 
\HH^1(G,\kk^{\times}) \rightarrow \HH^1(I,\kk^{\times})^{G/I} \rightarrow 
\HH^2(G/I,(\kk^{\times})^I) \,.
\end{equation}
Inertia acts trivially on $\kk$, and so the first and last groups in 
\eqref{eqJ}vanish by Hilbert's Theorem 90 and by the triviality of the 
Brauer group of finite fields, respectively.  Therefore
$$ \HH^1(G,\kk^{\times}) \cong \HH^1(I,\kk^{\times})^{G/I} \cong 
\Hom(I,\kk^{\times})^{G/I} \,.$$
The right side evidently has size no bigger than $e = \#I$, and so to 
complete our proof we need only show that the cocycles 
$g \mapsto \left(\frac{g\pi}{\pi}\right)^i$ for $0 \le i < e$ lie in 
distinct cohomology classes.  It suffices to show that they are 
nontrivial if $i \neq 0$.  However, if $b \in \kkux$ is such that 
$\left(\frac{g\pi}{\pi}\right)^i = \frac{b}{{}^g b}$ for all $g \in G$, 
then since the left-hand side of this equality has no terms involving $u$
we find $\left(\frac{g\pi}{\pi}\right)^i 
= \frac{b(0)}{g(b(0))}$ for all $g \in G$ as well.  Then $b(0) \pi^i \in 
K^G = L$, and by considering the valuation we see that $e \, | \, i$, so 
$i=0$.

\end{proof}

\section{Rank $1$ modules}

We retain our notation from the previous section, so in particular $K/L$
is a tamely ramified Galois extension of local fields with ramification
index $e = e(K/L)$, and $\pi \in K$ satisfies $\pi^e \in L$.
We let $e_K = e(K/\Ql)$ denote the
absolute ramification index of $K$, and $G =
\Gal(K/L)$ acts on $\kkuek$ as in Lemma \ref{galcoh}.  We will 
frequently need 
to divide various
integers by the greatest common divisor $(l-1,e)$, and so we make the
following defintion.
\begin{defn} {\rm
If $x$ is any integer, then $x'$ will denote $x/(l-1,e)$;  
moreover, use of the expression $x'$ will implicitly mean that $x$ is 
divisible by $(l-1,e)$.  
}
\end{defn}
It is often useful that $x'y=xy'$.
Finally, we choose $U$ an integer 
which is an inverse of $(l-1)'$ modulo $e'$, and let $V$ satisfy
\begin{equation}\label{eqK}
U(l-1)' = 1 + Ve' \,.
\end{equation}
When $l-1 \ | \ e$, for example, we will always choose $U=1$, $V=0$.

\begin{prop} \label{rankone} Consider the category of Breuil modules for 
$\OO_K$ with 
$\pi$ as the fixed choice of uniformizer.
\begin{enumerate}
\item[(1)]
A rank $1$ Breuil module $\M(r,a)$ admits generic fibre descent data from 
$K$ to $L$ 
if and only if $r$ is divisible by $(l-1,e)$ 
and $a \in {\mathbf l}^{\times} ((\kk[u]/u^{el})^{\times})^{l-1}$.

\item[(2)]
For each $0 \le r' \le e'_K$, $a \in {\mathbf 
l}^{\times}$, $c \in \Z/(l-1,e)\Z$, define $\M_U(r,a,c)$ to be the Breuil 
module $\M(r,a)$ with descent data given by $[g] \e = 
\left(\frac{g\pi}{\pi}\right)^{e' c - Ulr'} \e$, extended 
$G$-semilinearly to additive bijections on $\M$.  Any rank 
$1$ Breuil module with descent data from $K$ to $L$ is isomorphic to an 
$\M_U(r,a,c)$, and $\M_U(r,a,c) \cong \M_U(s,b,d)$ if and only if $r=s$, 
$c=d$, and $a/b \in ({\mathbf l}^{\times})^{l-1}$.
\end{enumerate}
\end{prop}

\begin{proof}
The first statement of the proposition follows immediately 
from the second, and from our understanding of maps between rank $1$ 
Breuil modules.  It is straightforward to check that the additive 
bijections in (2) do indeed define generic fibre descent data.

Suppose now that $\M(r,a)$ admits generic fibre descent data 
given by $[g]\e = A_g \e$.  We make the following observations:
\begin{enumerate}
\item We have $[h][g]\e = [h](A_g \e) = {}^h A_g A_h \e$, and so 
$A_{hg} = {}^h A_g A_h$.
\item Replacing $\e$ by $\etil = b\e$ as the standard basis vector, we 
find 
$\p(u^r \etil) = b^{l-1} a \etil$, and $[g]\etil = \frac{{}^g b}{b} A_g 
\etil$
\end{enumerate}
Therefore $g \mapsto A_g$ is a cocycle in $H^1(G,\kkuekx)$, while 
replacing $\e$ by $b \e$ multiplies this cocyle by the coboundary $g \mapsto 
\frac{{}^g b}{b}$.  As a consequence, by Lemma \ref{galcoh} we may make a 
choice of $\e$ so that $A_g = \left(\frac{g\pi}{\pi}\right)^{k}$ for 
some $k$.

Having done this, we calculate $$\p[g](u^r \e) = \p \left(
\left(\frac{g\pi}{\pi}\right)^{k+r} u^r \e \right) =
\left(\frac{g\pi}{\pi}\right)^{lk+lr} a \e$$ 
while
$$[g]\p(u^r \e) = [g] a\e = \left(\frac{g\pi}{\pi}\right)^{k} ({}^g a) \e 
\,.$$
Equating these two expressions, we find
$$ \left(\frac{g\pi}{\pi}\right)^{(l-1)k+lr} = \frac{{}^g a}{a} \,.$$
We conclude that $g \mapsto 
\left(\frac{g\pi}{\pi}\right)^{(l-1)k+lr}$ is a coboundary, and therefore 
$e \,| \, (l-1)k + lr$.  From this, we easily see that $(l-1,e) \,| \, r$ 
and therefore that $k$ must be of the 
form $e' c - U lr'$ for some $c \in \Z/(l-1,e)\Z$.

Since $e \,| \, (l-1)k+lr$ it also follows that ${}^g a/a = 1$ for all $g 
\in G$ and so $a \in (\kkuekx)^G = ({\mathbf l}[u^e]/u^{e_K l})^{\times}$.  
We 
note, however, that replacing $\e$ by $\etil = b\e$ with $b \in 
(\kkuekx)^G$ leaves the 
$A_g$ unchanged but replaces $a$ by $a b^{l-1}$.  Since 
$$(({\mathbf l}[u^{e}]/u^{e_K l})^{\times})^{l-1} = ({\mathbf 
l}^{\times})^{l-1}(1 + u^e {\mathbf l}[u^e]/u^{e_K l}) \,,$$ the remaining 
statements in the Proposition follow.
\end{proof}

\begin{example} {\rm  Recall our notation from Section \ref{dieudonne}.
We are interested in describing the rank $1$ Breuil
modules over $\OO_{E'}$ with descent data from $E'$ to $\Ql$ using a 
choice of $\pi = 
(-l)^{1/(l^2-1)}$ as our 
uniformizer.  In this case the ramification indices are $e = e_{E'} =l^2 - 
1$ and 
the residue 
field of~$E'$ is $\Fll$, so these are rank $1$ modules over
$\Fll[u]/u^{l(l^2-1)}$.  Since $l-1 \ | \ e$, we choose $U=1$ and $V=0$.
Whenever $U=1$, we drop the subscript $U$ from $\M_U(r,a,c)$.

Recall that
$\Gal(E'/\Ql)$ is 
generated by $\varphi$
and the $g_{\zeta}$, subject to the relations 
$g_{\zeta}^{l^2-1}=1$ and $\varphi
g_{\zeta} \varphi = g_{\zeta}^l$ for all $\zeta \in \Fll$.  
Since $E'/\Ql$ is tamely ramified and $\pi^e = -l \in \Ql$, 
we 
conclude from Proposition \ref{rankone} that the desired Breuil modules 
with descent data
are the $\M(r,a,c)$ given by
$$ \M = \langle \e \rangle , \ \ \MM = \langle u^{r'(l-1)} \e \rangle, \ \ 
\p(u^{r'(l-1)} \e) = a\e$$
for $a \in \Flx$ and $r=r'(l-1)$ with $0 \le r' \le l+1$, and 
generic fibre descent data 
given by $\fr(\e)=\e$ and 
$$[\zeta](\e) = \zeta^{(l+1)c - lr'} \e $$ with $c \in 
\Z/(l-1)\Z$.  (We will abbreviate $[g_{\zeta}]$ by $[\zeta]$.)  For 
different 
triples $(r,a,c)$, the $\M(r,a,c)$ are non-isomorphic.
}
\end{example}

\section{Identification of rank $1$ Breuil modules with descent data}

We now give an argument which identifies the Breuil modules with descent
data from the preceeding section, in the following sense: the Breuil
module
$\M_U(r,a,c)$ corresponds to a finite flat group scheme $\G_U(r,a,c)$
over $\OO_{K}$ with descent data, and we wish to determine the finite
flat group scheme $\G_{U,L}(r,a,c)$ over $L$ to which the generic fibre
$\G_U(r,a,c)_{/K}$
descends.  That is, 
we will compute the character $ \chi_{U,L}(r,a,c) : G_{L} \rightarrow
\Flx$ obtained as the Galois representation on $\G_{U,L}(r,a,c)(\Qlbar)$.
We will sometimes write $\chi_{U}(r,a,c)$ as shorthand for 
$\chi_{U,L}(r,a,c)$.

To begin, we note:

\begin{lemma} \label{maps} There is a non-zero homomorphism from 
$\M_U(r,a,c)$ to
$\M_U(s,b,d)$ if and only if $r \le s$ and $r \equiv s \pmod{l-1}$, $a/b 
\in ({\mathbf l}^{\times})^{l-1}$, and $c \equiv d + 
V(\frac{r-s}{l-1}) \pmod{(l-1,e)}$.
\end{lemma}

\noindent \textit{Remark:} Henceforth, we will always let $\e$ denote the 
standard basis vector of $\M_U(s,b,d)$ and $\ee$ the standard basis vector 
of $\M_U(r,a,c)$.

\begin{proof} Ignoring generic fibre descent data for the moment, the
descent-dataless analogue of this lemma (e.g. part 2 of Lemma 5.2.1 in
\cite{BCDT}) states that we have a non-zero map from $\M(r,a)$ to
$\M(s,b)$ if and only if $r \le s$, $r \equiv s \pmod{l-1}$, and $a/b \in 
(\kk^{\times})^{l-1}$, 
and moreover all such maps
are of the form $\ee \mapsto \alpha u^{l(\frac{s-r}{l-1})} \e$ for $\alpha 
\in \kk^{\times}$.
Such a map is compatible with generic fibre descent data exactly when
$$ \alpha \gpi^{e'c - Ulr'} u^{l(\frac{s-r}{l-1})} \e = g(\alpha) 
\gpi^{l(\frac{s-r}{l-1})} \gpi^{e'd - Uls'} u^{l(\frac{s-r}{l-1})} \e$$
for all $g \in G$.  This amounts to $\alpha \in {\mathbf 
l}^{\times}$ and 
$$ e'c - Ulr' \equiv e'd - Uls' + l\left(\frac{s-r}{l-1}\right) \pmod{e} 
\,, $$
and this congruence is easily seen to be equivalent to
$$ c - d \equiv V\left(\frac{r-s}{l-1}\right) \pmod{(l-1,e)}\,.$$
\end{proof}

\begin{cor} If $r \equiv s \pmod{l-1}$, then $\chi_U(r,a,c) = \chi_U(s,a,c 
+ V\left(\frac{s-r}{l-1}\right))$.
\end{cor}  

\begin{proof}  Put $d = c + V(\frac{s-r}{l-1}) \in 
\Z/(l-1,e)\Z$.  Suppose $r 
\le s$.  A non-zero 
map $\M_U(r,a,c) \rightarrow \M_U(s,a,d)$ exists by the previous lemma, 
and
corresponds to a non-zero map $\G_U(r,a,c) \rightarrow \G_U(s,a,d)$ 
compatible
with generic fibre descent data.  Therefore we get a non-zero map 
$\G_{U,L}(r,a,c) 
\rightarrow
\G_{U,L}(s,a,d)$.  This amounts to a non-zero map of Galois modules of
order $l$, so is therefore an isomorphism, and we find
$\chi_U(r,a,c)=\chi_U(s,a,d)$.  If $s > r$, the maps go in the other 
direction but the conclusion is the same.
\end{proof}

\begin{thm} \label{identify} Let 
$$\lambda = -la(\pi^e)^{Vr' - (l-1)'c} \in L \,.$$  
Then 
$\chi_U(r,a,c) : G_L \rightarrow \Flx$ is the character
$$ \eta_{\lambda} : g \mapsto 
\frac{g(\lambda^{1/(l-1)})}{\lambda^{1/(l-1)}}\,.$$
\end{thm}

Before we give the proof of Theorem \ref{identify}, we need the following 
lemma:

\begin{lemma} \label{scalar} Let $\G$ be a group scheme of order $l$ over 
$\OO_K$ with descent data from $K$ to $L$, such that by the Oort-Tate 
classification $\G \cong \Spec \OO_K[X]/(X^l + CX)$ with
comultiplication \eqref{eqF}.  
If $g \in \Gal(K/L)$, then the generic fibre descent data map $[g]$ 
sends $X \mapsto \alpha X$ for some $\alpha$ satisfying $\alpha^{l-1} = 
g(C)/C$.
\end{lemma}

\begin{proof} Let $\langle m \rangle$ be the multiplication-by-$m$ 
endomorphism of $\G$, and let $\chi : \Flx \rightarrow \Zlx$ denote the 
Teichm\"uller map.  Then $\langle m \rangle$ and $[g]$ commute, and so 
$[g]$ also commutes with the operator
$$ e_1 = \frac{1}{l-1} \sum_{m \in \Flx} \frac{1}{\chi(m)} \langle m 
\rangle$$
defined by Oort-Tate \cite{oorttate}.  But in the Oort-Tate 
construction, $e_1$ is the 
projection onto the submodule generated by $X$, and so
$$ [g](X) = [g] \circ e_1(X) = e_1 \circ [g](X) \in \OO_K \cdot X \,.$$
Since $[g](X^l) = \alpha^l X^l$ and $[g](CX)=g(C)\alpha X$, 
it follows that $\alpha^{l-1} = g(C)/C$.
\end{proof}

Now we return to the proof of Theorem \ref{identify}.

\begin{proof}  Let $K_0$ be the maximal unramified extension of $L$ 
inside $K$, and let $L_1 = L(\pi)$, so that $K=K_0 L_1$.  Then it 
suffices to show 
\begin{equation}\label{eqT}
\chi_{U,K_0}=\eta_{\lambda} \, |_{G_{K_0}}
\end{equation}
 and
and 
\begin{equation}\label{eqU}
\chi_{U,L_1} = \eta_{\lambda} \, |_{G_{L_1}}.
\end{equation}  
The formula \eqref{eqT} 
is precisely what is obtained by applying this theorem to the 
totally ramified 
extension $K/K_0$.  As for \eqref{eqU}, the extension $K/L_1$ is 
unramified, and the 
parameters $U_1$ and $V_1$ for this extension satisfy $U_1 (l-1) = 1 
+ V_1$.  Putting $$\lambda_1 = -la(\pi)^{V_1 r - (l-1)c} = -la(\pi)^{-r + 
(l-1)(rU_1 - c)}\,,$$ 
the statement of this theorem for the extension $K/L_1$ is that 
$\chi_{U,L_1}(r,a,c) = \eta_{\lambda_1}$.  But it is easily checked that 
$\lambda_1 / \lambda$ is a power of $\pi^{l-1}$, and so $\eta_{\lambda_1} 
= \eta_{\lambda} \, |_{G_{L_1}}$.  We conclude that it suffices to prove 
the theorem in the cases of $K/L$ an unramified extension and $K/L$ a totally 
(tamely) ramified extension.

The unramified case is easy, since in that case generic fibre descent data 
actually 
descends the 
group scheme.   Specifically, 
let $\M_L(r,a)$ be the Breuil module for $\OO_{L}$ with chosen 
uniformizer $\pi$ and parameters $r$ and $a$.  There is only one way 
to put generic fibre descent data on $\M(r,a)$, and 
Corollary 5.4.2 of \cite{BCDT} tells us that $\M(r,a)$ descends to 
$\M_L(r,a)$.  Then it follows from Lemma \ref{comult} that the affine 
algebra underlying $\M_L(r,a)$ is $\OO_L[X]/(X^l + \frac{la}{\pi^r}X)$ and 
$\chi_{U,L}(r,a) = \eta_{-la/\pi^r} = \eta_{\lambda}$, as desired.

We next turn to the situation when $L/K$ is totally (tamely) ramified.  We 
will first
consider the case when $e \, | \, r$.  Note that $\gpi^{e'} \in 
\Flx$, so 
that in this case $g$ acts as multiplication by 
$\gpi^{e'c - 
Ulr'} \in \Flx$ on the standard basis vector of $\M_U(r,a,c)$.  Since $K/L$ 
is totally ramified, the residue fields ${\mathbf l},\kk$ are equal, and 
part 2 of Theorem 5.6.1 of \cite{BCDT} tells us that $D([g])$ is the 
multiplication-by-$\gpi^{e'c - Ulr'}$ endomorphism 
of the Dieudonn\'{e}
module
$D(\G_U(r,a,c) \times_{K} {\mathbf k})$.  Since $\gpi^{e'c-Ulr'} \in \Fl$,
there is an integer $n$ such that $D([g]) = \Id + \cdots + \Id$ (where
there are $n$ $\Id$'s in the sum).  As the Dieudonn\'{e} functor is
additive, it follows that the corresponding action of $[g]$ on
$\G_U(r,a,c) \times_{K} {\mathbf k}$ is also multiplication by $n$,
and where $\gpi^{e'c-Ulr'}$ is the mod-$l$ reduction of $n$.

Recall from Lemma \ref{comult} that 
the
affine algebra of $\G_U(r,a,c)$ is $\OO_{K}[X]/(X^{l} + C X)$ where $C 
= \frac{\pi^{e_K-r} a}{G_{\pi}(\pi)} = \frac{la}{\pi^{r}}$,
and the comultiplication is
$$ X \mapsto 1 \otimes X + X \otimes 1 - \frac{l}{C} 
\sum_{i=1}^{l-1} \frac{X^i}{w_i} \otimes \frac{X^{l-i}}{w_{l-i}} \,.$$
In any case one can verify that the multiplication-by-$n$ map on the 
mod-$\pi$ reduction of this group scheme is $X \mapsto nX$, and using the 
fact that $n \equiv \gpi^{e'c-Ulr'} \pmod{\pi}$,
we therefore know that the action of $[g]$ on $\G_U(r,a,c)$ fits in 
the commutative diagram
$$ \xymatrix{
 {\G_U(r,a,c)} \ar[r]^{[g]} \ar[d]& {\G_U(r,a,c)} \ar[d]    \\ 
 {\Spec{\OO_{L}}} \ar[r]_{g} &   {\Spec{\OO_{L}}}  } $$
and sends $X \mapsto \gpi^{e'c-Ulr'} X + \pi f(X)$
for some polynomial $f(X)$.  Since $g(C) = C$, we conclude from Lemma 
\ref{scalar} that $[g] : X \mapsto \gpi^{e'c-Ulr'} X$.  

We now consider the action of this descent data on the generic fibre.
Put $\beta = \pi^{-(e'c-Ulr')}$ and $C_1 = \beta^{l-1} C \in L$.  
Note that we
have a $K$-algebra isomorphism $\gamma : K[X]/(X^{l} + C_1 X) 
\rightarrow
K[X]/(X^{l} +C X)$ sending $X \mapsto \beta X$,
and for each $g \in \Gal(K/L)$ we obtain a commutative diagram:
$$ \xymatrix{
 {K[X]/(X^{l} + C_1 X)} \ar[r]
\ar[d]_{X \mapsto X} & {K[X]/(X^{l} + C X)} \ar[d]^{[g]}  \\
 {K[X]/(X^{l} + C_1 X)} \ar[r] &
{K[X]/(X^{l} + C X)}  } $$
where the horizontal maps are $\gamma$ 
and the vertical maps are $g$-semilinear.
This shows that, pulling our generic fibre descent data back via the map 
$\gamma$,
the generic fibre descent data acts on $K[X]/(X^l + C_1 X)$ simply via
the action of Galois on $K$.  Moreover, $\lambda$ pulls back the 
comultiplication on
$K[X]/(X^{l} + C X)$ to the following comultiplication on
$K[X]/(X^{l} + C_1 X)$:
$$ X \mapsto 1 \otimes X + X \otimes 1 - \frac{l}{C_1}
\sum_{i=1}^{l-1} \frac{X^i}{w_i} \otimes \frac{X^{l-i}}{w_{l-i}} \,.$$
It follows immediately that the descended group scheme over $L$
corresponding to the above group scheme with descent data over $K$
is $L[X]/(X^{l} + C_1 X)$ with the usual Oort-Tate
comultiplication.  Therefore the character $\chi_{U,L}(r,a,c)$ is 
$\eta_{-C_1}$.  Using $e \, | \, r$ it is straightforward to check that 
$-C_1 / \lambda$ is an $(l-1)^{\textrm{st}}$ power in $L$, and so 
$\chi_{U,L}(r,a,c) = \eta_{\lambda}$ as well.

Now suppose instead that $0 < r < e_K$.  Again $\G_U(r,a,c)$ has
underlying algebra $\Spec \OO_K[X]/(X^l + CX)$ with $C =
\frac{la}{\pi^r}$.  From Lemma \ref{scalar} we have $[g] : X \mapsto \beta
X$ where $\beta$ satisfies $\beta^{l-1} = \gpi^{-r}$, and we must
determine $\beta$.

Note that the induced action of $[g]$ on the closed fibre sends $X \mapsto 
\beta X$, and by the identification of the Cartier-Manin Dieudonn\'e 
module of $\G_U(r,a,c) \otimes_{\OO_K} \kk \cong {\mathbf \alpha}_l$ with 
its tangent space, we find that $[g]$ induced multiplication by $\beta$ 
on the Dieudonn\'e module.  We also know that $[g]$ on 
$\M_U(r,a,c)/u\M_U(r,a,c)$ acts as 
multiplication by 
$\gpi^{e'c - lUr'}$, and from the proof of Theorem 5.1.3 in \cite{BCDT} it 
follows that $\beta^l = \gpi^{e'c - lUr'}$.  We conclude that $\beta = 
\gpi^{e'c - lUr' + r}$.  (One may check that indeed this $\beta$ satisfies 
$\beta^{l-1} = \gpi^{-r}$.) Now proceeding exactly as in the case $e \, | 
\, 
r$, we put 
$$C_1 = \pi^{-(e'C - lUr' + r)} \frac{la}{\pi^r} = la (\pi^e)^{(lVr' - 
C(l-1))'} \in L$$ 
and compute that $\chi_U(r,a,c) = \eta_{-C_1} = \eta_{\lambda}$.
\end{proof}

\begin{example} {\rm We return to the example of particular interest, 
namely 
when the extension $K/L$ is $E'/\Ql$, so $e_K = e = l^2-1$, $U=1$, 
$V=0$, and $\pi = (-l)^{1/(l^2-1)}$.  Then $(l-1)' = 1$, $\pi^e = -l$, and 
$\lambda = a(-l)^{1-c}$.  Theorem \ref{identify} now says that 
$\chi(r,a,c) = \chi_a \omega^{1-c}$, where $\chi_a$ is the unramified 
character sending Frobenius to $a$.}
\end{example}

\section{Rank $2$ extensions of rank $1$ modules} \label{ranktoo}

In this section, we classify the extensions, in the category of
Breuil modules with descent data, of the rank 
$1$ modules in the previous sections by one another.  The extensions 
without descent data are classified in Lemma 5.2.2 of \cite{BCDT}:

\begin{lemma} In the category of Breuil modules 
corresponding
to finite flat $l$-torsion group schemes over $\OO_{K}$ with choice of
uniformizer $\pi$, we have an isomorphism $$ \Ext^1
(\M(r,a),\M(s,b)) \cong
\{ h \in u^{\max(0,r+s-e_K)} \kkuek \} / \{ u^s t - (b/a) u^r t^l
\} $$
given by associating to each $h \in u^{\max(0,r+s-e_K)} \kkuek$ the 
$\p$-module
$$ \M = \la \e, \ee \ra, \ \MM = \la u^s \e, u^r \ee + h \e \ra, $$
with $$ \p(u^s \e) = b \e, \ \p(u^r \ee + h \e) = a \ee.$$  Moreover, 
replacing the basis element $\ee$ with $\ee + (b/a) t^l \e$ 
transforms $h$ to $h + (u^s t - (b/a) u^r t^l)$, and all equivalences 
between extensions are of this form.
\end{lemma}

We now wish to understand extensions of rank $1$ modules in the category of 
Breuil modules with descent data.  The underlying Breuil module extension
must be of the above form, and
 generic fibre descent data must act via
 $$ [g](\e)  = \gpi^{k_1} \e, \ [g](\ee) = \gpi^{k_2} \ee + A_g \e$$ 
where for ease of notation we have set
$$ k_1 = e'd - Uls' \ , k_2 = e'c - Ulr' \,.$$

One checks that the relation $[h][g]\ee = [hg]\ee$ is equivalent to
$$ \frac{A_{hg}}{(hg\pi/\pi)^{k_2}}
= \frac{A_{h}}{(h\pi/\pi)^{k_2}}
+ \left(\frac{h\pi}{\pi}\right)^{k_1-k_2} { }^h 
\left(\frac{A_{g}}{(g\pi/\pi)^{k_2}}\right)  $$
and so the map $g \mapsto \frac{A_{g}}{(g\pi/\pi)^{k_2}}$ is a cocycle in
$\HH^1(G,\kkuek)$ where $h \in G$ acts on $\kkuek$ via $h \cdot f = 
\left(\frac{h\pi}{\pi}\right)^{k_1-k_2} ({}^h f)$.  (The notation ${}^h f$ 
will always be reserved for the action defined in Lemma \ref{galcoh}.)  By 
the same method of proof as in Lemma \ref{galcoh}, this cohomology group 
vanishes, and so this map is in fact a coboundary.  

Now putting $\etil' = \ee + \frac{b}{a} t^l \e$, one computes 
\begin{eqnarray*}
[g] \etil' 
& = & \gpi^{k_2} \etil' + \left( A_g + \gpi^{k_1} 
{}^g \left(\frac{b}{a} t^l\right) - \gpi^{k_2} \frac{b}{a} t^l \right) 
\e \\
& = & \gpi^{k_2} \etil' + \left( A_g + \gpi^{k_2} \left( g \cdot \left( 
\frac{b}{a} t^l 
\right) - \frac{b}{a} t^l \right) \right)
\end{eqnarray*}
and so this alters $A_g/(g \pi/\pi)^{k_2}$ by 
the coboundary of
$\frac{b}{a} t^l$.  Since $A_g/(g \pi/\pi)^{k_2}$ is already a coboundary,
to see that in this fashion 
the $A_g$ may be transformed to $0$ by an appropriate choice of 
$t$ it suffices to show that all nonzero terms of $A_g$ have degree 
divisible by $l$.

To this end, we apply the relation $\p[g]=[g]\p$ to the element $u^r \ee + 
h\e \in \MM$.  One computes that
$$\p[g](u^r \ee + h\e) = \gpi^{(k_2+r)l} a \ee + 
\left(\frac{\Delta}{u^s}\right)^l b\e$$
where $$ \Delta = \gpi^r u^r A_g + {}^g h \gpi^{k_1} - h \gpi^{k_2 + r} $$
must have lowest term of degree at least $s$, while
$$[g]\p(u^r \ee + h\e) = \gpi^{k_2} a \ee + A_g a \e \,.$$
That the $\ee$-terms are equal follows from the fact that $e \ | \ 
(l-1)k_2 + lr$, while the equality between the $\e$-terms shows that $A_g$ 
is indeed an $l^{\textrm{th}}$ power.  

We can suppose, then, that all $A_g=0$.  Since now $(\Delta/u^s)^l$ 
must be $0$, 
it follows that a necessary and sufficient condition on $h \in 
u^{\max(0,r+s-e_K)} \kkuek$ for this 
extension of Breuil modules to admit generic fibre descent data is
$$u^{s+e_K} \ | \ \Delta = {}^g h \gpi^{k_1} - h \gpi^{k_2 + r} $$
for all $g \in G$.  Moreover, two such extensions with descent data 
with parameters $h,h'$ are equivalent precisely when $h'$ is 
of the form $h + u^s t - \frac{b}{a} u^r t^l$ for some $t$ such that 
$g \cdot \left(\frac{b}{a} t^l\right) = \frac{b}{a} t^l$ for all $g \in 
G$.  That is, we have the following necessary and sufficient conditions:
\begin{itemize}
\item all monomial terms of $h$ with degree $k < s + e_K$ must have
$k \equiv r + k_2 - k_1 \pmod{e}$ and coefficient in ${\mathbf l}$, and
\item all terms of degree $k < e_K$ of an allowable change-of-variables 
$t$ must have $k \equiv l^{-1} (k_2 - k_1)$ and coefficient in ${\mathbf 
l}$.
\end{itemize}

Before continuing, given $H \in \kkuek$ we describe an inductive 
procedure to 
solve the equation $H = u^s T - \frac{b}{a} u^r T^l$.  Let $H = \sum_i H_i 
u^i$ and $T = \sum_i T_i u^i$, so that the equation we wish to solve 
amounts to the 
system of equations 
\begin{equation}\label{eqL}
H_i = T_{i-s} - \frac{b}{a} T_{\frac{i-r}{l}}^{l}
\end{equation}
for $0 \le i < l e_K$, and where $T_j$ is required to be $0$ if $j$ is 
not a nonnegative integer.  Set $i_0 = \frac{ls-r}{l-1}$.  We will 
attempt to solve 
the equations \eqref{eqL} inductively, inducting on the distance of 
$i$ from $i_0$.  The condition $|i - i_0| < \frac{1}{l-1}$ is an empty 
condition unless $i = i_0$ is a nonnegative integer, in which case the 
associated equation is
$$ H_{\frac{ls-r}{l-1}} = T_{\frac{s-r}{l-1}} - \frac{b}{a} 
T_{\frac{s-r}{l-1}}^l \,.$$
If this equation can be solved for $T_{\frac{s-r}{l-1}}$, this is our base 
case, and then assume the following inductive hypothesis:
\begin{itemize}
\item  the equations \eqref{eqL} 
can be solved for all $i$ such that $|i - i_0| < 
N + \frac{1}{l-1}$, and
\item  in doing so, all and only the $T_j$ with $|j - 
\left(\frac{s-r}{l-1}\right)| < N + \frac{1}{l-1}$ have been determined.
\end{itemize}

Now suppose that $i$ satisfies
$N + \frac{1}{l-1} \le |i - i_0| < N + \frac{l}{l-1}$.
Then $N + \frac{1}{l-1} \le |(i-s) - \left(\frac{s-r}{l-1}\right)| < N + 
\frac{l}{l-1}$ and so by assumption $T_{i-s}$ has not been determined.
On the other hand, $|\frac{i-r}{l} - \frac{s-r}{l-1}| < \frac{N}{l} + 
\frac{1}{l-1} \le N + \frac{1}{l-1}$, and so $T_{\frac{i-r}{l}}$ 
\textit{has} been determined.  So we may recursively take
$$T_{i-s} = H_i + \frac{b}{a} T_{\frac{i-r}{l}}^{l} \,.$$
This is only a condition if $i < s$, in which case there is a solution 
only if the $T_{i-s}$ so-obtained is $0$.  By induction, we conclude that 
the system \eqref{eqL}
has a solution if and only if:
\begin{itemize}
\item the base case $H_{\frac{ls-r}{l-1}} = T_{\frac{s-r}{l-1}} - 
\frac{b}{a} T_{\frac{s-r}{l-1}}^l $ is either vacuous or is non-vacuous 
and has a solution, and
\item in our recursive process, $T_{i-s} = H_i + \frac{b}{a} 
T_{\frac{i-r}{l}}^{l} = 0$ for whenever $i < s$.
\end{itemize}  Note that the base case may be unsolvable only if 
$\frac{s-r}{l-1}$ is a negative integer and $H_{\frac{ls-r}{l-1}} \neq 0$;
or if $\frac{a}{b} \in (\kk^{\times})^{l-1}$ so that the map $\alpha 
\mapsto 
\alpha - \frac{b}{a} \alpha^l$ is not surjective.  In the latter case, 
fix any $\eta$ not in the image of the map $\alpha \mapsto \alpha - 
\frac{b}{a} \alpha^l$.   

As an example of the usefulness of this description, we can employ it
show:

\begin{prop} \label{uniq} Suppose $H = u^s T - \frac{b}{a} u^r T^l$ has a 
solution and 
$\deg H < s$.  Then $H=0$.
\end{prop}

\begin{proof}  If the base 
case is not vacuous, then either $\frac{ls-r}{l-1} \ge s$ and so 
$H_{\frac{ls-r}{l-1}} = 0$ 
by assumption, or else $r > s$ and the 
assumption that the equation can be solved forces
$H_{\frac{ls-r}{l-1}} = 0$; in any case the base case may be solved 
by taking $T_{\frac{s-r}{l-1}} = 0$.    We claim that in our inductive 
procedure, all $T_i$ will be determined to be $0$: indeed, if $i \ge s$,
then $T_{i-s} = H_i + \frac{b}{a} T_{\frac{i-r}{l}}^l = 0$ by induction, 
while if $i < s$ then perforce $T_{i-s} = 0$. Thus if the system of 
equations \eqref{eqL} can be 
solved, then $T=0$ is a solution, and so $H=0$.
\end{proof}

In a similar vein, we can show
\begin{prop} \label{killer} Let $H$ be as before.

\begin{enumerate}
\item If the base case for $H$ is vacuous, can be solved, or 
cannot be solved but $r > s$,  then there exists a unique $H'$ such that 
$H' = u^s T - \frac{b}{a} u^r 
T^l$ can be solved and such that $\deg(H-H') < s$.  

\item If the 
base case cannot be solved and $s \ge r$, then there exists a unique $H'$ 
such that the only nonzero term of $H-H'$ of degree at least $s$ is of the 
form $N \eta u^{\frac{ls-r}{l-1}}$ for $N \in \Fl$.

\end{enumerate}

\end{prop}

\begin{proof} For part (a), uniqueness is evident by Proposition \ref{uniq}.  
Existence when the base case is vacuous or can be solved follows from the 
inductive procedure for $H$, simply defining $H'_i = -\frac{b}{a} 
T_{\frac{i-r}{l}}^l$ whenever $i<s$.  If $r > s$ and the 
base case cannot be solved, first set $H'_{\frac{ls-r}{l-1}} = 0$, and 
proceed as before.  For part (b), if $s > r$ and the base case cannot be 
solved, then since the $N\eta$ are coset representatives for $\{ \alpha - 
\frac{b}{a} \alpha^l \}$ in $\kk$  there is a unique 
$H'_{\frac{ls-r}{l-1}} \in \{ \alpha - \frac{b}{a} \alpha^l \}$ such that 
$H_{\frac{ls-r}{l-1}} - H'_{\frac{ls-r}{l-1}} = N\eta$ for some $n$, and 
then we proceed to construct $H'$ via the inductive procedure as in part 
(a).
\end{proof}

Finally we return to the situation under consideration, namely that 
all monomial terms of $h$ with degree $k < s + e_K$ must have
$k \equiv r + k_2 - k_1 \pmod{e}$ and coefficient in ${\mathbf l}$.
Notice that if $i - s 
< e_K$, then $i < s + e_K$ and $\frac{i-r}{l} < e_K$.  Moreover, since 
$k_2 - k_1 \equiv s - r + l^{-1}(k_2 - k_1) \pmod{e}$, we find that 
$i - s \equiv l^{-1}(k_2 - k_1) \pmod{e}$ if and only if $\frac{i-r}{l} 
\equiv l^{-1}(k_2-k_1) \pmod{e}$ if and only if $i \equiv k_2 - k_1 + r 
\pmod{e}$.

We use our procedure to attempt to solve the equation $h 
= u^s t - 
\frac{b}{a} u^r t^l$.  
Suppose first that the base case is vacuous, 
or cannot be 
solved but $r >
s$.   Using the above observations, and by induction, 
the coefficient $t_{i-s}$ for 
$i - s < e_K$ can become 
nonzero only if $i \equiv s + l^{-1}(k_2-k_1) \pmod{e}$; and in 
that 
case 
induction and the formula for $t_{i-s}$ in terms of $H_i$ and 
$t_{\frac{i-r}{l}}$ shows that $t_{i-s} \in {\mathbf l}$.  It 
follows from the method of part (a) of Proposition \ref{killer}
that when we construct $h'$ and $t$ such that $h'= u^s t - 
\frac{b}{a} u^r t^l$ and $\deg(h-h')<s$, the resulting $t$ satisfies our 
condition that every nonzero term of degree $k$ smaller than $e_K$ has $k 
\equiv 
l^{-1} (k_2 - k_1)$ and coefficient in ${\mathbf l}$.  Moreover, also by 
construction, the terms of $h'$ of degree less than $s + e_K$ have 
coefficients in ${\mathbf l}$, and so all the coefficients of $h-h'$ lie 
in ${\mathbf l}$.

Next, consider the situation where $s \ge r$, and the base case is 
non-vacuous, so $s \equiv r \pmod{l-1}$.  If 
$\frac{ls-r}{l-1} \not\equiv r + k_2 - k_1 \pmod{e}$, then 
$h_{\frac{ls-r}{l-1}} = 0$ and so taking $t_{\frac{s-r}{l-1}}=0$ the 
conclusions of the previous paragraph hold.  Suppose, then, that 
$\frac{ls-r}{l-1} \equiv r + k_2 - k_1 \pmod{e}$, or in other words 
that $\frac{s-r}{l-1} \equiv l^{-1}(k_2 - k_1) \pmod{e}$. It is not 
difficult 
to see that this congruence is equivalent to $c-d \equiv 
V\left(\frac{r-s}{l-1}\right) \pmod{(l-1,e)}$.   If $a/b$
is not an $(l-1)^{\textrm{st}}$ power in ${\mathbf l}^{\times}$, then 
the base case can solved with $t_{\frac{s-r}{l-1}} \in 
{\mathbf l}$, and again the conclusions of the previous paragraph hold.  
(Note that the congruence $\frac{s-r}{l-1} \equiv l^{-1}(k_2 - k_1) 
\pmod{e}$ ensures that the possibly-nonzero coefficient 
$t_{\frac{s-r}{l-1}}$ lies in suitable degree.)

We are finally left with the case when $s \ge r$, $s \equiv r \pmod{l-1}$,
$a/b$ is an $(l-1)^{\textrm{st}}$ power in ${\mathbf l}^{\times}$, and $c 
-d \equiv V\left(\frac{r-s}{l-1}\right) \pmod{(l-1,e)}$.  Note that this 
is exactly the case when there is 
a nontrivial map $\M_U(r,a,c) \rightarrow \M_U(s,b,d)$.  
Let $\eta$ be any fixed
element of ${\mathbf l}$ not in the image of $\alpha \mapsto \alpha - 
\frac{b}{a} \alpha^l$.  Following the 
method of part (b) of Proposition \ref{killer} and using the same 
arguments as in the previous paragraphs, we construct $t$ such that 
the nonzero terms of $t$ of degree $k < e_K$ have $k \equiv l^{-1} (k_2 - 
k_1)$ and coefficient in ${\mathbf l}$, and such that $h - (u^s t - 
\frac{b}{a} t^l)$ has coefficients in ${\mathbf l}$ and all terms of 
degree 
less than $s$, save possibly for a term of the form $N\eta 
u^{\frac{ls-r}{l-1}}$.  Putting this all together, we obtain:

\begin{thm} \label{ranktwo} Put $k_1 = e'd - Uls'$ and $k_2 = e'c- Ulr'$.
\begin{enumerate}
\item Suppose that there is no map $\M_U(r,a,c) \rightarrow \M_U(s,b,d)$.  
Then every extension of $\M_U(r,a,c)$ by $\M_U(s,b,d)$ with descent data 
is isomorphic to exactly one of the form:
$$\M = \langle \e,\ee \rangle , \ \ \MM=\langle u^s\e, u^r \ee + h\e 
\rangle \,,$$
$$ \p(u^s \e) = b\e , \p(u^r \ee + h\e) = a\ee \,,$$
$$ [g](\e) = \gpi^{k_1} \e, [g](\ee) = \gpi^{k_2} \ee \,,$$
where $h \in u^{\max(0,r+s-e_K)} {\mathbf l}[u]/u^{e_K l}$ has 
degree less than $s$ and all nonzero terms of degree congruent to $r + k_2 
- k_1 \pmod{e}$.  In particular, the $\Fl$-dimension of this space of 
extensions is at most $[L:\Ql]$.
\item If there is a map $\M_U(r,a,c) \rightarrow \M_U(s,b,d)$, let $\eta$ 
be any fixed element of ${\mathbf l}$ not in the image of $\alpha 
\mapsto \alpha - \frac{b}{a} \alpha^l$ on ${\mathbf l}$.  Then the same 
conclusion holds as in part (a), except that $h$ may also have term of the 
form $N\eta u^{\frac{ls-r}{l-1}}$.  In particular, the $\Fl$-dimension of 
this space of extensions is at most $[L:\Ql]+1$.
\end{enumerate}
\end{thm}

To see the dimension claims, note in part (a) that at most 
$e_K/e = e_L$ 
different terms in $h$ can be nonzero.  Since each coefficient lies in 
${\mathbf l}$, the dimension over $\Fl$ is at most $e_L f_L = [L:\Ql]$.  
The claim in part (b) follows identically.

We remark that this result is intuitive: the number of extensions grows as 
$s$ gets larger and $r$ gets smaller, in other words as the group scheme 
corresponding to $\M_U(s,b,d)$ gets ``more \'etale'' and that 
corresponding to $\M_U(r,a,c)$ gets ``more multiplicative''.  This is 
sensible as there are plenty of extensions of \'etale group schemes by 
multiplicative ones, and none in the other direction.

\subsection{Start of the proof of Theorem \ref{rest}} \label{start}

We return once again to the case when $K/L$ is the extension $E'/\Ql$.
Since we are only interested in $2$-dimensional representations of
$G_{\Ql}$ with nontrivial centralizer, we are safely in the situation
where there is no map $\M(r,a,c) \rightarrow \M(s,b,d)$, for otherwise the
two diagonal characters would be equal; that is, since $V = 0$ we are
assuming $(a,c) \neq (b,d)$.  In this case, part (a) of Theorem
\ref{ranktwo} tells us that the space of extensions with descent data
of $\M(r,a,c)$ by
$\M(s,b,d)$ is at most $1$-dimensional, and in fact nonsplit extensions
exist exactly whenever there is a solution to the congruence 
$$ n \equiv (l+1)(c-d) + ls' - r'$$ 
with 
$$ \max(0,r+s-(l^2-1)) \le n < s \,.$$ 
Then write $h = h_n u^n$, and denote the resulting extension
$\M(r,a,c;s,b,d;n,h_n)$.

Since we have assumed that $\rhobar$ has centralizer $\kk$,
we may henceforth restrict ourselves to the nonsplit situation above;
in particular $s \neq 0$ and $r \neq l^2-1$, so that an $n$
satisfying the given inequality and congruence can exist.  
Moreover, we
will always take $h_n = 1$: if $h_n
\neq 1$, the resulting group scheme with descent data
 is isomorphic to the group scheme with descent data having
identical parameters save $h_n = 1$ --- they are simply non-isomorphic as
extension classes, which will be of no concern.
Therefore, we will need to consider only Breuil modules with descent data
 of the form
$\M(r,a,c;s,b,d;n,1)$.

We now turn to the question of which of the
group schemes with descent 
data corresponding to these Breuil modules with descent data
satisfy the
relations (\ref{eqC}, \ref{eqD}, \ref{eqE}) 
on their Dieudonn\'{e} modules: namely, that $\z +
\z^l$ acts as $\zeta^{-m} + \zeta^{-lm}$, that $\z^{l+1}$ acts as
$\zeta^{-(l+1)m}$, and that $F + TV = F + abV$ acts at $0$.  
It is easy to see,
using the compatibility between Dieudonn\'{e} theory and
Breuil theory described in Section \ref{bm}, that
the Dieudonn\'{e} module of the closed fibre of $\M(r,a,c;s,b,d;n,1)$ has
a basis $\mathbf{v},\mathbf{w}$ on which $\z$ acts in the
following manner:
$$\z(\mathbf{v}) = \zeta^{(l+1)c-lr'} \mathbf{v} \ \text{and} \ 
  \z(\mathbf{w}) = \zeta^{(l+1)d-ls'} \mathbf{w}$$
and so
$$\z^l(\mathbf{v}) = \zeta^{(l+1)c-r'} \mathbf{v} \ \text{and} \
  \z^l(\mathbf{w}) = \zeta^{(l+1)d-s'} \mathbf{w}.$$ 
It follows that if $\z$ satisfies the desired relations, then either
$$  \zeta^{(l+1)c-lr'}=\zeta^{-m} \ \text{and} \
\zeta^{(l+1)c-r'}=\zeta^{-lm}$$
or
$$ \zeta^{(l+1)c-lr'}=\zeta^{-lm} \ \text{and} \ 
\zeta^{(l+1)c-r'}=\zeta^{-m},$$ and a similar relationship holds between
$d$, $s'$, and $m$.  Recalling that $m = i + (l+1)j$, the first
possibility
yields congruences
$$(l+1)c - lr' \equiv -m \pmod{l^2-1} \ \text{and} \
  (l+1)c - r' \equiv -lm \pmod{l^2-1}.$$
Solving for $r'$, we obtain $ r' \equiv -i \pmod{l+1} $.  Since $1 \le i
\le l$ and $0 \le r' \le l+1$, we conclude $r' = l + 1 - i$.  This allows
us to solve that $c \equiv 1 - i - j \pmod{l-1}$ (which completely 
determines $c$, since it is an element of $\Z/(l-1)\Z$).  Applying a
similar analysis to the second of the possible sets of relations among
$c$, $r'$, and $m$, we find that between the two cases,
$$ (r',c) = (i, -j) \ \text{or} \ (l+1-i,1-i-j) .$$
By an identical calculation,
$$ (s',d) = (i, -j) \ \text{or} \ (l+1-i,1-i-j) .$$
However, if $r'=s'$ and $c=d$, we would require 
$$n \equiv (l+1)(c-d) + ls' - r' \equiv (l-1)s' = s \pmod{l^2-1}.$$
Since we require $0 \le n < s$ and since $s < l^2-1$ (as $i \neq l+1$),
this situation is impossible.  We have therefore proved that the only
possibilities for Breuil modules with descent data
attached to integral 
models with descent data for
$\rhobar$ which satisfy our Dieudonn\'{e} module relations
are those of the form
$$ \M((l-1)(l+1-i),a,1-i-j; (l-1)i, b, -j; 0,1) $$ and
$$ \M((l-1)i,a,-j; (l-1)(l+1-i), b, 1-i-j; 0,1). $$
By Theorem \ref{identify} (and recalling the contravariance of
the Breuil module functor) the descended $G_{\Ql}$-representations
$\rhobar$ corresponding to these Breuil modules with descent data
 are exactly of the form
\begin{itemize}
\item
$\rhobar =
\begin{pmatrix}
 \omega^{i+j} \chi_a & * \\
    0         & \omega^{1+j} \chi_b
\end{pmatrix}
$ and
\item
$\rhobar =
\begin{pmatrix}
 \omega^{1+j} \chi_a & * \\
    0         & \omega^{i+j} \chi_b
\end{pmatrix} 
$.
\end{itemize}
Notice that unless $i=1$ or $i=l$, each different possibility for
$\rhobar$ yields at most one Breuil module with descent data
 in our list.  When
$i=1$ and $i=l$, it is still possible that
$\M(l(l-1),a,-j;(l-1),b,-j;0,1)$ and $\M((l-1),a,-j;l(l-1),b,-j;0,1)$
are both integral models with descent data
for the same $\rhobar$; however, we will see in the next
section that the former arises from a residual representation of $G_{\Ql}$
which is either split or is nonsplit but does not have centralizer $\Fl$, and 
since we
have assumed that $\rhobar$ has
centralizer $\Fl$ this group scheme with descent data
 cannot arise from our $\rhobar$.
Therefore it is again the case that our $\rhobar$ gives rise to
at most one integral model with descent data.  We will also prove in the
next section that if $\rhobar =
\begin{pmatrix}
 \omega^{i+j} \chi_a & * \\
    0         & \omega^{1+j} \chi_b
\end{pmatrix}
$ gives rise to one of the integral models with descent
data in the above list and if $i=2$, then $*$ is
peu-ramifi\'{e}.  Similarly if $\rhobar =
\begin{pmatrix}
 \omega^{1+j} \chi_a & * \\
    0         & \omega^{i+j} \chi_b
\end{pmatrix}$
and if $i=l-1$, then $*$ is peu-ramifi\'{e}.  Once done, all of these
results together will have
completed the proof of
\begin{prop} \label{rightlist} If $\tau = \omt^{m} \oplus \omt^{lm}$ and
$\rhobar : G_{\Ql}
\rightarrow \GL(\Fl)$ has centralizer $\Fl$ and is reducible, and if
$\Rdvz \neq 0$, then $\rhobar \, |_{I_l}$ does indeed have one of the
forms specified in Theorem \ref{main}.  Furthermore $\rhobar$ gives rise
to exactly
one finite flat group scheme over $\OO_{E'}$ with descent data to $\Ql$
satisfying the necessary relations on the Dieudonn\'{e} module of its
closed fibre.
\end{prop}

\begin{remark}{\rm The relation $F + T V = F + abV = 0$ is indeed
satisfied on
the Dieudonn\'{e} modules of the closed fibres of the above group
schemes.  One may check that on our basis $\mathbf{v}$,$\mathbf{w}$, $F$
and $V$ act via the matrices
$$\begin{pmatrix}
-b & 0 \\
 0 & 0
\end{pmatrix} \ \text{and} \
\begin{pmatrix}
1/a & 0 \\
 0  & 0
\end{pmatrix}$$ respectively.}
\end{remark}

\section{Maps between rank $2$ Breuil modules with descent data}

\subsection{Generalities}

Our strategy for proving that certain pairs of rank $2$ Breuil modules
with descent data arise from the same representation is to find maps
between these rank $2$ modules.  

\begin{defn} {\rm  Let $\M$ be the Breuil module corresponding to a 
group scheme $\G$ over $\OO_K$.  Then, by Raynaud \cite{Raynaud},
$\G$ is mapped to by a maximal integral model $\G_{+}$ and maps to
a minimal integral module $\G_{-}$.  The maximal and minimal Breuil module
of $\M$ are defined to be, respectively, the Breuil modules corresponding
to $\G_{+}$ and $\G_{-}$.
}\end{defn}

By Lemma 4.1.4 of \cite{BCDT}, if two
extensions with descent data 
of rank $1$ Breuil modules for $\OO_K$ with descent data from 
$K$ to $L$
arise from the same representation of $G_{L}$,
they both map to a maximal Breuil module with descent data for this
representation, and are also mapped to by a minimal Breuil module with descent
data, where the maps are generic fibre isomorphisms.  
A scheme-theoretic
closure argument as in Lemma 4.1.3 of \cite{BCDT} shows that in this 
situation, the maximal
and minimal Breuil module with descent data
are again extensions of rank $1$ Breuil modules with
descent data.  

As a first example, note that as a corollary of \ref{maps} and
\ref{identify}, if $\chi_U(r,a,c)=\chi_U(r',a',c')$ then there is a
nonzero map either from $\M_U(r,a,c)$ to $\M_U(r',a',c')$ or vice-versa,
depending on whether $r \le r'$ or $r' \le r$.  For instance, when $l-1 \,
| \, e$, so that $U=1$ and $V=0$, this implies that the maximal and
minimal Breuil modules with descent data
of $\M(r,a,c)$ are $\M(e_K,a,c)$ and $\M(0,a,c)$
respectively.  It follows easily that if $\M$ is an extension with descent
data of $\M(r,a,c)$ by $\M(s,b,d)$, where $\chi(r,a,c) \neq \chi(s,b,d)$,
then the descended $G_L$-representation of $\M$ is split if and only if
there is a nonzero map $\M(0,a,c)  \rightarrow \M$, and that in this case 
the 
maximal and minimal Breuil modules with descent data
 are $\M(e_K,a,c)\oplus\M(e_K,b,d)$ and 
$\M(0,a,c)\oplus\M(0,b,d)$ respectively.  In a similar vein:

\begin{prop} \label{mapsexist} Suppose that we have a diagram 
$$ 
\xymatrix{ 
\G'_1 \ar[r] & \G_1 \ar[r] \ar[d] & \G''_2 \\ 
\G'_2 \ar[r] & \G_2 \ar[r] & \G''_2 }
$$ 
where for $i=1,2$, $\G'_i$ and $\G''_i$ are finite flat group schemes 
over $\OO_K$ of order $l$ with descent data from $K$ to $L$ whose
generic fibres descend to non-isomorphic irreducible $G_L$-representations, 
and where $\G_i$ 
is an extension with compatible generic fibre descent data.
  Suppose 
furthermore 
that the map $\G_1 \rightarrow \G_2$ induces a generic fibre isomorphism
of group schemes with descent data,
and that the descended generic fibre representation of 
$G_L$ is not semisimple.  Then there are maps 
$\G'_1 \rightarrow \G'_2$ and 
$\G''_1 \rightarrow \G''_2$ which are isomorphisms on the generic fibre.
\end{prop}

\begin{proof}  By the semisimplicity assumption, the 
irreducible $G_L$-representations
 corresponding to $\G'_1$ 
and $\G''_2$ are different.  Therefore the 
composite map $$\G'_1 \rightarrow \G_1 
\rightarrow \G_2 \rightarrow \G''_2$$ is the zero map and so factors 
through $\G'_2$.  Let $\check{\G}$ denote the Cartier dual of $\G$.  
Dualizing 
our diagram, we obtain a nonzero map $\check{\G}''_2 
\rightarrow \check{\G}''_1$, and dualizing again gives a nonzero map 
$\G''_1 
\rightarrow \G''_2$.
\end{proof}

Finally, we will need to make use of the following result:

\begin{prop} \label{genericisom} Let $\G \rightarrow \G'$ be a map of 
finite flat
group schemes over $\OO_K$ of
equal order, both killed by $l$.
If the kernel of the corresponding map $\M' \rightarrow \M$ of Breuil
modules does not
contain a free $\kk[u]/u^{el}$-submodule, then $\G \rightarrow \G'$ is
an isomorphism on generic fibres.
\end{prop}

\begin{proof}  Assume $f: \G \rightarrow \G'$ does not induce an isomorphism
on generic fibres.  
Then the image of $f_{/K}$ in $G'_{/K}$
is not all of $G'_{/K}$, and taking scheme-theoretic closure of this image 
yields an exact sequence of group schemes
$$ 0 \longrightarrow {\mathcal H} \longrightarrow \G' \overset{g}{\longrightarrow} {\mathcal H}' \longrightarrow 0 $$
with ${\mathcal H}' \neq 0$ and $g \circ f = 0$.  
If $\N'$ is the Breuil module corresponding to ${\mathcal H}'$,
then $\N' \hookrightarrow \ker(\M' \rightarrow \M)$, since short-exact
sequences of group schemes yield short-exact sequences of Breuil modules.
\end{proof}

Note that if the map $\G \rightarrow \G'$ in the preceeding proposition is 
in fact a map of group schemes with descent data, then the isomorphism in 
the conclusion is also an isomorphism of group schemes with descent data.  

\subsection{Application to the proof of Proposition \ref{rightlist}}
\label{rktwo}

Return once again to the situation where the extension $K/L$ is $E'/\Ql$, 
and suppose henceforth that $(a,c) \neq (b,d)$, so that our 
representations have different diagonal characters.  We begin with the 
following:

\begin{prop} The descended representation of 
$\M(r,a,c;s,b,d;n,1)$ is split if and only if $r > s$ and $c=d$.
\end{prop}

\begin{proof} By the discussion in the previous section, we must determine 
when there exists a nonzero map $$ \Psi: \M(0,a,c) \rightarrow 
\M(r,a,c;s,b,d;n,1) \,.$$  Let $\f$ denote the standard basis vector of 
$\M(0,a,c)$.  Assume that $\Psi$ exists. Then $\Psi(\f) = V \e + W \ee$, 
and since $(a,c)\neq(b,d)$ it follows that $W \neq 0$.  

From the fact that $\Psi$ commutes with generic fibre descent data, it 
follows that all
nonzero terms of $V$ are in degrees congruent to $(l+1)(c-d) + ls'
\pmod{l^2-1}$ and have coefficients in $\Fl$; and that all nonzero terms of
$W$ are in degrees congruent to $lr' \pmod{l^2-1}$ with coefficient in
$\Fl$.  From the fact that $\Psi$ commutes with $\p$, it follows that all
nonzero terms of $V$ and $W$ are in degrees divisible by $l$.  By the
Chinese remainder theorem, $V = vu^{\alpha}$ and $W = w u^{\beta}$ are
monomials with $v,w \in \Fl$, and from the given conditions it follows 
with one exception that $\alpha = 
l((l+1)\{c-d\} + s')$ and $\beta=lr'$, where for $x \in \Z/(l-1)\Z$, 
$\{x\}$ is the unique representative of $x$ lying between $0$ and $l-2$.
The exception is that $\alpha =0$ when $s'=l+1$ and $\{c-d\}=l-2$, as in 
this case $l((l+1)\{c-d\} + s')=l(l^2-1)$.

Since $V\e + W\ee \in \MM(r,a,c;s,b,d;n,1)$, it follows that 
\begin{equation} \label{eqM}
\beta \ge r \,, \ \alpha \ge s \,, \ \beta - r + n \ge s \,,
\end{equation}
 and 
$$V \e + W\ee = (vu^{\alpha-s} - wu^{\beta-r+n-s}) u^s\e + 
wu^{\beta-r}(u^r \ee + u^n \e) \,.$$  Note that the inequalities \eqref{eqM}
rule 
out the possibility $s'=l+1$ and $\alpha=0$, and so we indeed have
$\alpha=l((l+1)\{c-d\} + s')$.  The condition that $\Psi$ commutes with 
$\p$ is then equivalent to
$$ bv u^{l(\alpha-s)} - bw u^{l(\beta-r+n-s)} = av u^{\alpha} \,.$$

Now $\beta - r + n - s = r' + n - s < r'$, so $l(\beta-r+n-s) < 
l(l^2-1)$ and the term $bw u^{l(\beta-r+n-s)}$ is nonzero.  It follows 
that $v$ is nonzero, and since a sum of three monomials can equal zero 
only if each nonzero term in the sum has the same degree, we see that 
$\alpha = l(\beta-r+n-s)$.  This yields $n = (l+1)\{c-d\} + ls' - r'$, and 
since $n < s$ this forces $c=d$ and $r' > s'$.

Finally, under the assumptions $c=d$ and $r'>s'$, one can check 
that
$$ \f \mapsto vu^{ls'} \e + \left(1 - \frac{a}{b}\right) vu^{lr'} \ee$$
is a map of the desired sort, and so these conditions are sufficient as 
well as necessary.
\end{proof}

Note that as claimed in Section \ref{start}, this shows that 
descended $G_{\Ql}$-representation corresponding to
the Breuil 
module $\M(l(l-1),a,-j;(l-1),b,-j;0,1)$ is split.

Observe that the preceeding proposition may be reinterpreted as follows.  
The quantity $(l+1)\{c-d\} + ls' - r'$ lies between $0$ and $2(l^2-1)$,
and so $n = (l+1)\{c-d\} + ls' - r' - N(l^2-1)$ for some $N \in 
\{0,1,2\}$.  The equality $n = (l+1)\{c-d\} + ls' - r'$ occurs precisely 
when $\{c-d\}=0$ and $r > s$, and this is exactly the case when the 
descended $G_{\Ql}$-representation is split.  
Note that $(l+1)\{c-d\} + ls'-r' = 
2(l^2-1)$ only when $\{c-d\}=l-2$, $s=l^2-1$, and $r=0$, and the rest 
of the time we must have $n = (l+1)\{c-d\} + ls' - r' - (l^2-1)$.  We can 
now prove:

\begin{prop} \label{ranktwomaps}
Suppose that $$\M(r,a,c;s,b,d;n,1) \ \textrm{and}  \ 
\M(r_1,a_1,c_1;s_1,b_1,d_1;n_1,1)$$ have non-split descended 
$G_{\Ql}$-representation.  Then there is a map
$$ \M(r,a,c;s,b,d;n,1) \rightarrow \M(r_1,a_1,c_1;s_1,b_1,d_1;n_1,1) $$
which is an isomorphism on generic fibers if and only if 
$(a,c)=(a_1,c_1)$, $(b,d)=(b_1,d_1)$, $r \le r_1$, and $s \le s_1$.
\end{prop}

\begin{proof} The conditions in the proposition are necessary by 
Proposition \ref{mapsexist}.  To see that the conditions are sufficient, 
we will exhibit the desired maps
$$\M(r,a,c;s,b,d;n,1) \rightarrow \M(r_1,a,c;s_1,b,d;n_1,1)$$
whenever $r \le r_1$ and $s \le s_1$.  Let $\e,\ee$ and $\f,\ff$ denote 
our standard bases for the left-hand and right-hand Breuil modules 
with descent data
above, respectively.  Most of the time, we have
$$ n - (ls' -r') = n_1 - (ls'_1 - r'_1) = (l+1)\{c-d\} - (l^2-1) $$
and in these cases one can check that the maps given by 
\begin{eqnarray*}
\e & \mapsto & vu^{l(s'_1-s')} \f \\
\ee & \mapsto & vu^{l(r'_1-r')} \ff
\end{eqnarray*}
are indeed maps of Breuil modules with descent data.  The equality 
$n-(ls'-r')=n_1-(ls'_1-r'_1)$ is crucial to the verification that the map
preserves the filtration and commutes with $\p$.  Similarly, we have maps
$$ \M(0,a,d-1;s,b,d;ls'-(l+1),1) \rightarrow \M(0,a,d-1;l^2-1,b,d;0,1)$$
given by
\begin{eqnarray*}
\e & \mapsto & vu^{l(l+1-s')} \f \\
\ee & \mapsto & a^{-1}bv \f + a^{-1}bv \ff
\end{eqnarray*}
and maps
$$ \M(0,a,d-1;l^2-1,b,d;0,1) \rightarrow 
\M(r_1,a,d-1;l^2-1,b,d;l^2-1-r'_1,1)$$
given by
\begin{eqnarray*}
\e & \mapsto & v \f \\
\ee & \mapsto & - v \f + b^{-1}av u^{lr'_1} \ff \,.
\end{eqnarray*}
This exhibits the desired maps in the remaining cases, when one or the 
other Breuil module with descent data has $\{c-d\}=l-2$, $s=l^2-1$, and 
$r=0$.

To see that these maps all induce isomorphisms on the generic fibre, we 
note that this follows from the following general criterion.  If we have a 
map of such Breuil 
modules sending 
\begin{eqnarray*}
\e & \mapsto & v u^{\alpha} \f \\
\ee & \mapsto & y u^{\beta} \f + z u^{\gamma} \ff
\end{eqnarray*}
with $v,z \neq 0$ 
then every 
element of the kernel of our homomorphism is annihilated by 
$u^{\max(\gamma, \alpha + \gamma - \beta)}$. 
If $\alpha + \gamma - \beta < l(l^2-1)$, this shows that the kernel does 
not contain any free $\Fll[u]/u^{l(l^2-1)}$-submodules, and so by 
Proposition \ref{genericisom} the map induces an isomorphism on generic 
fibres.
\end{proof}

The analogous result is true for maps between Breuil modules with descent
data in which the descended $G_{\Ql}$-representation is split.

\subsection{Lattices of rank $2$ Breuil modules}

Consider the non-split Breuil module with descent data
$\M=\M(r,a,c;s,b,d;n,1)$.  Implicit 
in the existence of this Breuil module with descent data 
is that $n$ satisfies the 
inequalities $\max(0,r+s-(l^2-1)) \le n < s$ as well as the congruence $n 
\equiv (l+1)(c-d) + ls' - r' \pmod{l^2-1}$.  By our earlier comparisons 
of $n$ and $(l+1)\{c-d\} + ls' - r'$, we see that $n - 
(ls'-r') = -k(l+1)$ for an integer $k$ between $0$ and $l$.  Indeed, if 
$c=d$ then 
$k=0$ or $l-1$, the former if and only if $r > s$; if $c=d-1$, then $k=1$ 
or $l$, the latter if and only if $s=l^2-1$ and $r=0$; and otherwise $k$ 
is the unique integer between $2$ and $l-2$ which represents $d-c$.  In 
fact, we can show:

\begin{prop} \label{classify} For fixed $k$, the pairs $(r',s')$ 
for which $n = (ls'-r') - k(l+1)$ satisfies the inequalities 
$\max(0,r+s-(l^2-1)) \le n < s$ are 
precisely
the pairs satisfying:
$$ 0 \le r' \le l-k $$
and $$ k+1 \le s' \le l+1 \, $$
with the exceptions that for $k=0$ we require $r' > s'$, and for $k=1$
the pair $(0,l+1)$ is excluded. 
\end{prop}

\begin{proof}  We shall prove that for fixed $k > 0$ the desired
pairs $(r',s')$ 
are the lattice points inside the 
convex quadrilateral bounded by the
inequalities
$$ r' \ge 0 $$
$$ s' \le l+1 $$
$$ ls' - r' \ge (l+1) k$$
$$ s' - lr' + (l^2 - 1) \ge (l+1) k \,,$$
and, if $k = 1$, the region excludes the extremal point $(0,l+1)$.
For $k = 0$, we shall similarly prove that the pairs $(r',s')$ for 
which $n=(ls'-r')-k(l+1)$
are precisely the lattice points inside the triangle
bounded by the inequalities 
$$ r' > s' $$
$$ ls' - r' \ge (l+1) k$$
$$ s' - lr' + (l^2 - 1) \ge (l+1) k \,.$$
It is easy to see that the lattice points inside these regions are
exactly the ones described in the statement of the proposition.

 At the outset, we know that we must satisfy the inequalities 
$$0 \le s' \le l+1 \,,$$
$$0 \le r' \le l+1 \,,$$
$$ 0 \le n < s \,,$$
$$ r + s - (l^2-1) \le n \,.$$
From $n = (ls'-r')-k(l+1) < s$ we get $-k(l+1) < r'-s'$.  This is no
condition if $k \ge 2$, but if $k=1$ we exclude $(0,l+1)$ and if $k=0$
we need $r'>s'$.  The condition $n \ge 0$ translates into $ls'-r' \ge (l+1)k$, and the condition $n \ge r+s-(l^2-1)$ translates into
$s-lr' + (l^2-1) \ge (l+1)k$.  Therefore, the conditions in the
proposition are necessary.  We need to show that they are sufficient.

If $k=0$, the inequalities $ls' \ge r' > s'$ imply $s > 0$, so $r > 0$,
while the inequalities $r' + (l^2-1) > s' + (l^2-1) \ge lr'$ imply
$r' < l+1$, so $s' < l+1$.

If $k>0$, the inequalities $r' \ge 0$ and $ls' - r' \ge (l+1)k > 0$
imply $s > 0$, while the inequalities $s' \le l+1$ and $s' - lr'
+ (l^2-1) \ge (l+1)k > 0$ imply $r' < l+1$.
\end{proof}

\begin{figure}[ht]
\begin{center}
\epsfig{file=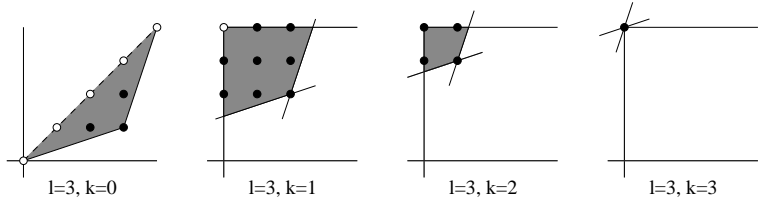,bbllx=0pt,bblly=0pt,bburx=407pt,bbury=105pt,%
        width=10cm}
\end{center}
\caption{The regions of Proposition \ref{classify} when $l=3$}
\end{figure}

\begin{cor} For all choices of $r$,$s$, and $n$
such that $$\M(r,a,d-1;s,b,d;n,1)$$ is a Breuil module with descent data
whose 
descended $G_{\Ql}$-representation is nonsplit, that representation is 
peu-ramifi\'{e}.
\end{cor}

\begin{proof} This is the $k=1,l$ case in the previous proposition.
The above discussion, combined with the maps constructed in 
Proposition \ref{ranktwomaps}, shows that the module 
$\M(r,a,d-1;s,b,d;n,1)$ has
minimal Breuil module with descent data
 $\M(0,a,d-1;2(l-1),b,d;l-1,1)$, and so for fixed $a,b,d$ 
these Breuil modules with descent data
all correspond to integral models with descent
data having the same descended
representation.  To see that
this representation is peu-ramifi\'{e}, we note (see, e.g., Section 8 
of \cite{edix}) that peu-ramifi\'{e} representations of $G_{\Ql}$ have
integral $\Zl$-models.  Therefore at least one of the above Breuil 
modules with descent data
corresonds to an integral model with descent data for a peu-ramifi\'{e}
representation.  Consequently, they all do.
\end{proof}

This completes the proof of Proposition \ref{rightlist}.  We summarize 
these results as follows:

\begin{thm} Fix $a,b \in \Flx$ and $c,d \in \Z/(l-1)\Z$, and let $\rho$
be a representation
$$\begin{pmatrix}
\chi_a \omega^{1-c} & * \\ 
0 & \chi_b \omega^{1-d}
\end{pmatrix} $$ of $G_{\Ql}$, with $* \neq 0$ and $(a,c) \neq (b,d)$.  
If $d-c \equiv 1 \pmod{l-1}$, suppose $*$ is peu-ramifi\'e.
Let $k$
be the integer between $1$ and $l-1$ congruent to $d-c \pmod {l-1}$.
Then the Breuil modules with descent data corresponding to the 
integral models with descent data for $\rho$ over $\OO_{E'}$ are
the Breuil modules with descent data 
$\M(r,a,c;s,b,d;n,1)$ with $$ 0 \le r' \le l-k $$
and $$ k+1 \le s' \le l+1.$$  
The lattice of these Breuil modules with descent data
is a square with $l-k+1$ points on
each side, and maps 
from  $\M(r,a,c;s,b,d;n,1)$ to  $\M(r',a,c;s',b,d;n',1)$
respecting generic fibre descent data exist whenever
$r' \ge r$ and $s' \ge s$.  In particular, there are $(l-k+1)^2$ such
integral models with descent data, and the maximal and minimal 
integral models for this 
representation correspond to the Breuil modules with descent data
$$ \M( (l - k)(l-1), a,c ; (l+1)(l-1), b, d; l^2 - kl ,1) $$
and
$$ \M( 0, a, c; (k+1)(l-1), b, d; l-k, 1) \,. $$
If $*$ were tr\`es-ramifi\'e, then $\rho$ would have no such integral
models with descent data.
\end{thm}

\begin{proof}  
Our analysis in Section \ref{rktwo} shows that the descended 
$G_{\Ql}$-representations 
of these Breuil modules with descent data
are indeed non-split.  The rest of the 
claims follow from Propositions \ref{ranktwomaps} and \ref{classify}.
\end{proof}

\section{Rank $4$ calculations}

Recall that our list of nonsplit rank $2$ Breuil modules with descent data
satisfying
Dieudonn\'{e} module conditions (\ref{eqC}, \ref{eqD}, \ref{eqE}) was
$$ \M((l-1)(l+1-i),a,1-i-j; (l-1)i, b, -j; 0,1) $$ and
$$ \M((l-1)i,a,-j; (l-1)(l+1-i), b, 1-i-j; 0,1). $$
Notice that the change of variables $i \mapsto l+1-i$ and $j \mapsto 
i+j-1$ interchanges the two collections of Breuil modules with descent data
above, so it 
suffices to consider the latter; moreover, we need consider only 
those Breuil modules with descent data
whose descended $G_{\Ql}$-representation is non-split and 
has nontrivial centralizer.  So, to prove Theorem \ref{rest} we are
reduced to showing, for each 
\begin{equation}\label{eqV}
\M = \M((l-1)i,a,-j; (l-1)(l+1-i), b, 1-i-j; 0,1)
\end{equation}
with $i = 1,\ldots,l-1$ and $j \in \Z/(l-1)\Z$, and $a \neq b$ if $i=1$, 
that the space of
extensions of $\M$ by $\M$ with descent data
still satisfying the desired Dieudonn\'{e}
module relations is at most $1$-dimensional.  We now begin this
computation. For clarity we will continue to write $r'$ for $i$ and $s'$
for $l+1-i$, since that is what we are used to.  Note that $r'+s'=l+1$,
$ls'-r'=(l-i)(l+1)$, and $lr'-s' = (l+1)(i-1)$.

Let $(\N,\NN,\p)$ be an arbitrary extension of $\M$ by $\M$ with descent 
data.  We will let
$\e,\ee$ denote the standard basis for the submodule $\M$ of $\N$, while
$\f,\ff$ will denote lifts of the standard basis for the quotient $\M$ of
$\N$.  Then $$ \N = \langle \e,\ee,\f,\ff \rangle, $$ and a priori $\NN$ 
has the form
$$ \NN = \langle u^{s} \e, u^{r} \ee + \e,  u^{s} \f +
A \e + B \ee, u^{r} \ff + \f + C \e + D \ee \rangle \,.$$
First, we wish to see that $\f,\ff$ may be chosen appropriately so that 
$A=B=C=D=0$.  To begin, replace $u^{s} \f +  A \e + B \ee$ 
with  $  u^{s} \f + A \e + B \ee -
A( u^{r} \ee  + \e)$ in our basis for $\NN$, so that we may take $A=0$.  
Similarly we can take $C=0$.   Now note that since $u^{l^2-1} \f \in
u^{l^2-1} \N \subset \NN$ and $u^{r}(u^{s} \f + B
\ee) \in \NN$, we obtain $u^{r} B \ee \in \NN$.  This implies
$u^{s} \, | \, B$.  Writing $B = u^{s} B'$, we may take
$\tilde{\f} = \f + B' \ee$, and then $\NN$ has the basis $\langle u^s \e, 
u^r \ee + \e, u^s \tilde{\f} , u^r \ff + f + D \ee \rangle$ for some $D$.
Finally, we wish to alter $\ff$ to eliminate $D$.  By the same
considerations as before, we can see $u^{s} D \ee \in \NN$, and so
$u^{r} \, | \, D$.  Putting $D = u^{r} D'$ we can take
$\tilde{\f}' = \ff + D' \ee$, and we conclude that we may suppose 
$$ \NN = \langle u^{s} \e, u^{r} \ee + \e,  u^{s} \f , 
u^{r} \ff + \f  \rangle \,.$$
The next thing we want to do is determine the ways we can still alter 
$\f$, $\ff$ to $\tilde{\f}$, $\tilde{\f}'$ 
while preserving this form for $\NN$, i.e., keeping
$u^{s} \tilde{\f}$, $u^{r} \tilde{\f}' + \tilde{\f} \in \NN$.  
To this end, suppose 
$$ \f \mapsto \tilde{\f} = \f + A \e + B' \ee \,, \ \ff \mapsto
\tilde{\f}' = \ff + C \e + D' \ee \,.$$
Then
$$ u^{s} \tilde{\f} = u^{s} \f + A u^{s} \e + B'
u^{s} \ee, $$ and this is in $\NN$ provided $u^{r}$ divides
$B'$.  Write $B' = u^{r} B$.  Now
$$ u^{r} \tilde{\f}' + \tilde{\f} = (u^{r} \ff + \f) + (A +
u^{r} C - B - D') \e + (B + D')( u^{r} \ee + \e) \,.$$
Thus $C$ may be arbitrary so long as we select $D'$ such that
$u^{s}$ divides $A + u^{r} C - B - D'$.  Writing 
$$ A + u^{r} C - B - D' = u^{s} D $$ we may evidently make
$D$ arbitrary and put $$D' =  A + u^{r} C - B -  u^{s} D
\,.$$  So our most general change of variables is 
\begin{equation}\label{eqN}
\tilde{\f} = \f + A \e + u^{r} B \ee \,, \ \tilde{\f}' = \ff + C
\e + (A - B +  u^{r} C +  u^{s} D) \ee
\end{equation}
with $A,B,C,D$ arbitrary.
We now turn to the question of $\p$.  We suppose
$$ \p(u^{s} \e) = b \e \,, \p(u^{r} \ee + \e) = a \ee $$
$$ \p(u^{s} \f) = b \f + V\e + W\ee \,, \p(u^{r} \ff + \f) = a
\ff + Y\e + Z\ee .$$  
Using the change-of-variables 
\eqref{eqN}, we wish to simplify $V$,$W$,$Y$,$Z$.  To begin with, we
try  $\tilde{\f} = \f + A \e + u^{r} B \ee$ (with a commensurate
choice of $\tilde{\f}'$, which for now will be irrelevant).
Then one computes that $\p(u^{(l-1)s'} \tilde{\f})$ is equal to
$$ b \tilde{\f} + (V - b A + b(A-B)^l) \e + 
 (W +  u^{r} ( a B^l u^{(ls-r)} - b B)) \ee \,.$$
Since $B$ may be arbitrary and $ls - r > 0$ we may make $ ( a B^l
u^{(ls-r)} - b B)$ arbitrary, and we may use this choice to
eliminate all terms in $W$ of degree at least $r$.  Thus we may assume
$\deg(W) < r$.  Making this change completely determines $u^{r} B$,
so we may now make this change and assume henceforth that
$B = 0$ and $\deg(W) < r$.  Then $V$ is altered to $V +
b(A^l - A)$ by our choice of $A$, which we can use to eliminate every term
of $V$ except the constant term.  We can therefore suppose that $V$
is a constant $v$, that $W$ is a polynomial of degree less than $r$, and 
that the only still-allowable change of $\f$ is $\f \mapsto \f + 
\alpha \e$, with
$\alpha$ a constant, moving $V \mapsto V + b(\alpha^l - \alpha)$.

Consider the additive map from $\Fll \rightarrow \Fll$ sending $x$ to $x^l
- x$.  The kernel is exactly $\Fl$, while if $x^l - x \in \Fl$ then $(x^l
- x)^l = x^l - x$; since $x^{l^2}=x$ and $l \neq 2$ we find $x^l = x$.
So our map induces an isomorphism $\Fll/\Fl \rightarrow
\Fll/\Fl$.  Thus we may select $\alpha$ above so that $V \in \Fl$, and
then $V$ is completely fixed, while $\f \mapsto \f + \alpha \e$ with
$\alpha \in \Fl$ is the only possible change of $\f$.

To reduce further, we now wish to see the ways in which these
extensions of Breuil modules admit generic fibre descent data.
Suppose $$[g] \f = \gpi^{(l+1)d - ls'} \f + A_{g} \e +
B_{g} \e' .$$
Then $$[g](u^s \f) = \left(\frac{g\pi}{\pi} u\right)^s  \left(\gpi^{(l+1)d 
- ls'} \f  + A_{g} \e
+ B_{g} \e' \right) \in \NN$$ which requires $u^r \, | \, B_{g}$, say
$B_{g}
= u^r B'_{g}$.  We see that $\p([g](u^s \f))$ is equal to
$$ \gpi^{(l+1)d - ls'} (b\f + v \e + W \ee) + \gpi^{ls} u^{sl}
(B'_{g})^l a \ee + g \gpi^{sl} (A_{g} - B'_{g})^l \e $$
whereas $[g](\p(u^s \f))$ is
$$ \gpi^{(l+1)d-ls'} b \f + (v \gpi^{(l+1)d-ls'} + b A_{g}) \e +
({}^g W \gpi^{(l+1)c-lr'} + b B_{g}) \ee $$
using ${}^g v = v$ since $v \in \Fl$.  Matching coefficients we get
\begin{equation}\label{eqR}
A_{g} = \gpi^{sl} (A_{g} - B'_{g})^l
\end{equation}
and
\begin{equation}\label{eqO}
{}^g W \gpi^{(l+1)c-lr'} + b B_{g} = w \gpi^{(l+1)d - ls'} + a 
\gpi^{sl} u^{sl} (B'_{g})^l \,.
\end{equation}

Since $W$ is of degree
less than $r$ whereas $B_{g}$ and $u^{sl}$ are divisible by $u^r$, 
\eqref{eqO} implies
\begin{equation}\label{eqP} 
{}^g W \gpi^{(l+1)c-lr'} =  W \gpi^{(l+1)d - ls'}
\end{equation}
and
\begin{equation}\label{eqQ} 
b B_{g} = a \gpi^{sl} u^{sl} (B'_{g})^l \,.
\end{equation}
All of the rank $2$ Breuil modules under consideration here have $n=0$,  
so $(l+1)d - ls' \equiv (l+1)c - r' \pmod{l^2-1}$, and using this in 
\eqref{eqP}
 we obtain ${}^g W = \gpi^{r} W$ for all $g$.  Since 
$\deg{W} < r$, this is only possible if $W=0$.  In \eqref{eqQ}, if 
the left-hand side has
lowest nonzero term of degree $k$ , then for the right-hand side the
lowest term has degree $sl + l(k-r)$.  Equating these degrees gives 
$k = r + (r'-ls') < r$, contradicting our divisibility condition on
$B_{g}$.  Thus $B_{g}=0$.  Taking $B'_{g}=0$ in \eqref{eqR}, we finally 
obtain
$$ A_{g} = \gpi^{sl} A_{g}^l $$
which implies that $A_{g}$ is a constant.  Indeed 
$A_{g} \gpi^{-(l+1)d + ls'} \in \Fl$, and one checks from this and from 
the relation $[h][g]=[hg]$ that the
map $g \mapsto \gpi^{-(l+1)d + ls'} A_{g}$ is a 
homomorphism from $G=\Gal(E'/\Ql)$ to $\Fl$, which must be the zero map.   
We have thus shown
$A_{g}=B_{g}=0$ and $W=0$.

Next we consider the more difficult problem of  
simplifying $Y$,$Z$, and $[g] \ff$ by altering $\ff$.  Taking
$A=B=0$ in \eqref{eqN}, we select
$$ \tilde{\f}' = \ff + C \e + (u^{r} C + u^{s} D) \ee \,.$$
Then one computes that $\p ( u^{r} \tilde{\f}' + \f )$ is
equal to
$$ a \tilde{\f}' + (Y - a C - b D^l) \e + (Z - a (u^{r} C + u^{s} D) 
+ a (u^{r} C + u^{s} D)^l ) \ee \,.$$
Whatever $D$ is, we will certainly want to take $a C = Y - bD^l$ to
eliminate $Y$ (which completely determines $C$ in terms of $D$).  We may
therefore assume $Y=0$ and $a C = - bD^l$, and then our map alters
$$ Z \mapsto 
Z  - ( b u^{r} D^{l}  - a u^{s} D)^l   + (b  u^{r} D^{l} - a u^{s}
D) \,.$$
Noting this, we now turn to the consideration of generic fibre
descent data.  Suppose
for each $g$ that 
$$[g] \ff = \gpi^{(l+1)c-lr'} \ff + E_{g} \e + F_{g} \ee \,.$$
We then
have 
$$[g](u^r \ff + \f) = \gpi^{(l+1)c-r'}(u^r \ff + \f) + \gpi^r (E_{g} 
u^r - F_{g}) \e + \gpi^r F_{g} (u^r \ee + \e)$$
so $u^s \, |  \, \gpi^r (E_{g} u^r -  F_{g}) = u^s \Delta_{g}$
and
$$\p [g] (u^r \ff + \f) =  \gpi^{(l+1)c-lr'} (a \ff + Z \e') +
\Delta_{g}^l b \e + a \gpi^{lr} F_{g}^l \ee \,.$$
Matching coefficients with
$$ [g] \p (u^r \ff + \f) = \gpi^{(l+1)c-lr'} a \ff + 
a E_{g} \e + a F_{g} \ee + {}^g Z \gpi^{(l+1)c-lr'} \ee $$ gives
$$ a E_{g} = b \Delta_{g}^{l} $$ and 
\begin{equation}\label{eqS} a F_{g} +
{}^g Z \gpi^{(l+1)c-lr'} = a \gpi^{lr} F_{g}^l +
Z \gpi^{(l+1)c-lr'} \,.
\end{equation}

With these equations in hand we compute from $[hg]=[h][g]$ that the map $g 
\mapsto 
\gpi^{-(l+1)c + lr'} E_g$ is a cocycle in the group cohomology 
$\HH^1(G,\coeffFll)$ where $G$ acts on $\coeffFll$ via $g \cdot f = \gpi^r 
{}^g f$.  Similarly $g \mapsto \gpi^{-(l+1)c + lr'} F_g$ is a cocycle in 
$\HH^1(G,\coeffFll)$ for the action $g \cdot f = {}^g f$.  We know from 
the 
proof of Lemma \ref{galcoh} that both these cohomology groups are trivial,
and therefore we obtain elements $P, Q \in \coeffFll$ such that
$$  \gpi^{-(l+1)c+lr'} E_{g} = \gpi^r {}^g Q - Q $$
and
$$  \gpi^{-(l+1)c+lr'} F_{g} =  {}^g P - P \,.$$
Setting $R = u^r Q - P$, we see
$$ E_{g} u^r - F_{g} = \gpi^{(l+1)c-lr'} ( {}^g R - R ) \,.$$
Recalling that $u^s \Delta_{g} = \gpi^r (E_{g} u^r - F_{g})$,
so that $$ \Delta_{g} = \gpi^{(l+1)c-r'} \frac{{}^g R - R}{u^s} $$ 
we find that $R$ must have no terms of degree less than $s$, except
possibily for a constant term in $\Fl$. 
We write $R = r_0 + r_s u^s + \cdots = r_0 + u^s R_0$.  Then the
equation $a E_{g} = b \Delta_{g}^l$ gives us:
$$  \left(\frac{{}^g R - R}{u^s}\right)^l = \frac{a}{b} 
\left(\gpi^r {}^g Q - Q\right) \,.$$
Writing $Q = q_0 + q_1 u + \cdots$, we examine the above equation
term-by-term.  Using that $r \equiv ls \pmod{l^2-1}$, the left-hand side
has terms of the form $({}^g r_{i+s}^l \gpi^{il + r} - r_{i+s}^l) u^{il}$  
while 
the
right-hand side has terms of the form $\frac{a}{b} ({}^g q_j \gpi^{r+j} - 
q_j) u^j$.  Thus $q_j = 0$ unless $j$ is divisible by $l$, or unless $j
\equiv s \pmod{l^2-1}$ and $q_j \in \Fl$.  If $j = il$ and is not 
congruent to $s
\pmod{l^2-1}$ then the map $x \mapsto {}^g x \gpi^{r+il} - x$ is injective 
and we can match $q_{il} = \frac{b}{a} r_{i+s}^l$.  
From this analysis, we conclude that
$$ Q = \frac{b}{a} R_0^l + Q' $$
where the terms of $Q'$ have degree congruent to $s \pmod{l^2-1}$ and 
coefficients in $\Fl$.  Therefore
$$  P  =  u^r Q - R = \frac{b}{a} u^r R_0^l - u^s R_0 + P'$$
where $P'$ has terms of degree divisible by $l^2-1$ and coefficients in 
$\Fl$.
Combining \eqref{eqS} with $$F_{g} = \gpi^{(l+1)c-lr'} ( {}^g P - P)$$
we get $$ {}^g (a P - a P^l + Z) = a P - a P^l + Z $$ and so 
all terms of $aP - aP^l + Z$
are of degree
divisible by $l^2- 1$ with coefficient in $\Fl$.
Putting all this together, we find that
$$Z =  (b u^r R_0(u)^l - a u^s R_0(u))^l - ( b u^r
R_0(u)^l - a u^s R_0(u)) + Z'$$
where $Z'$ has terms of degree divisible by $l^2-1$ and coefficients in 
$\Fl$.  Therefore \emph{taking $D = R_0$ in our change-of-variables \eqref{eqN}
for 
$\ff$ transforms
$Z$ into $Z'$}, a polynomial with all terms of degree divisible by 
$l^2-1$ and coefficients in $\Fl$.

We still wish to reduce $Z$ further, which is easier now that we can
assume $Z$ has no terms of low degree except a constant term in $\Fl$.  If 
we alter
$Z$ via some choice of $D$, we suppose that $D = \sum_i d_i u^i$ has no
terms of degree less than $s'-r'$.  Then the lowest nonzero term of
$$ \frac{b}{a} u^{lr} D^{l^2} - (\frac{b}{a} u^{r} + u^{ls}) D^l +  u^{s}
D$$
has degree $ls'-r'$, and specifically the lowest term is
$$ (-\frac{b}{a} d_{s'-r'}^l + d_{s'-r'}) u^{ls'-r'} \,. $$
The equation $x = (-\frac{b}{a} d_{s'-r'}^l + d_{s'-r'})$ may
be solved for 
$$  d_{s'-r'} = \left(1 -
\frac{b^2}{a^2}\right)^{-1} \left( x + \frac{b}{a}
x^l\right) $$ except possibly if $a = \pm b$ and $x \neq 0$.  Note
also that if $x \in \Fl$, there is a solution for $d_{s'-r'}$ 
except 
possibly if $a=b$.

The terms of degree $i > ls'-r'$ in our transformation of $Z$ are
$$ u^{s} (d_{i-s} u^{i-s}) - u^{ls} (d_{\frac{i-ls}{l}} 
u^{\frac{i-ls}{l}})^l  - \frac{b}{a} u^r (d_{\frac{i-r}{l}} 
u^{\frac{i-r}{l}})^l
+ \frac{b}{a} u^{lr} (d_{\frac{i-lr}{l^2}} u^{\frac{i-lr}{l^2}})^{l^2} $$

Since $i - s > (i-ls)/l$, $(i-r)/l$, $(i-lr)/l^2$ for $i > ls'-r'$
we see that taking $d_{i-s} = 0$ for $i$ up to $ls'-r'$ and solving
the resulting {\em linear} equations for $d_{i-s}$ for $i > ls'-r'$,
we may alter $Z$ to remove all terms of degree
greater than $ls'-r'$ (without introducing a term of degree $ls'-r'$ if
there wasn't one to begin with).  Therefore unless $ls'-r'=l^2-1$,
i.e. unless $r'=1$,$s'=l$, we may certainly take $Z$ to be a constant.
In case
$r'=1$,$s'=l$, note that the case $a=b$ is excluded automatically from our 
list of rank $2$ Breuil modules with descent data \eqref{eqV}, 
and so again the term of degree $ls'-r'=l^2-1$ can 
be removed by this argument.  Therefore in any case we can suppose $Z$ is 
a constant $z \in \Fl$.

In case $a = \pm b$ and $s' \ge r'$, let $\eta$ be a choice of
$(a/b)^{1/(l-1)}$.  
We note, for future reference, that for and $d \in \Fl$ by the above
argument there is a change-of-$\ff$ leaving $Z$ fixed,
given by $D = \eta d u^{s'-r'} + \text{(higher terms)}$ and the 
corresponding $C$.

Now observe that because we have reduced $Z$ to a simple form, we get 
$F_g = a \gpi^{lr} F_{g}^l$, and since $u$ divides $F_{g}$
we get $F_{g}=0$.  Then our equation for $E_{g}$ becomes 
$$ a E_{g} = b \left(\frac{\gpi^r u^r E_{g}}{u^s}\right)^l $$ and 
so if $E_{g}$ is nonzero then: $E_g$ is a monomial 
of degree $u^{l(s'-r')}$, but also $s' \ge r'$ and $a/b$ is an $(l-1)^{\rm 
st}$ power in $\Fllx$, i.e. $a = \pm b$.  

So automatically $E_g = 0$ 
unless $s' \ge r'$ and $a = \pm b$, which is exactly the situation in 
which there was a change-of-$\ff$ leaving $Z$ fixed.  In this case write 
$\gpi^{-(l+1)c + lr'} E_g = e_g u^{l(s'-r')}$, so that 
$e_{hg} = e_h + \left(\frac{h\pi}{\pi}\right)^{ls'-r'} {}^h e_g$.  If 
$s'=l,r'=1$ then $g \mapsto e_g$ is a homomorphism $G \rightarrow \Fll$, 
so is zero.   Otherwise, by the 
usual cohomology argument $e_g = \gpi^{ls'-r'} {}^g e - e$ for some $e \in 
\Fll$.  Then selecting any $g$ for which $\gpi^{ls'-r'} \in \Fl$ is not 
$1$, we see that $e = \eta d \in \eta \Fl$, where $\eta$ was our 
previously-chosen $(l-1)^{\rm st}$ root of $a/b$.  Finally, we make the 
change-of-$\ff$ which fixes $Z$ and has $D = \eta d u^{s'-r'} +$ (higher 
terms).  The corresponding $C = - a^{-1} b D^l = - \eta d
u^{l(s'-r')} +$ (higher terms).  Then $\tilde{\f}' = \ff + C \e + (\ee \ 
\text{term})$, and
we compute that $[g] \tilde{\f}'$ is equal to
$$\gpi^{(l+1)c-lr'} \tilde{\f}' + E_{g} \e +\left({}^g C 
\gpi^{(l+1)c-r'} - C \gpi^{(l+1)c-lr'} \right)\e +
(\ee \ \text{term}) \,.$$
So we have transformed $E_{g}$ into
$$ E_g - \left({}^g \eta d \gpi^{l(s'-r')} \gpi^{(l+1)c-r'} - \eta d 
\gpi^{(l+1)c-lr'}\right) u^{l(s'-r')}  + ({\rm higher}) \,,$$
and therefore this transformation leaves the $E_g$
with no term of degree $u^{l(s'-r')}$.
However, notice that since $Z$ is unchanged by this transformation
we still obtain $F_{g} = 0$, and now our new
$E_{g}$, having no terms of degree $u^{l(s'-r')}$, is also $0$.

To summarize, we have proved

\begin{thm}  For the $\M$ under consideration, any $\N \in \Ext^1 (\M,\M)$
in the category of Breuil modules with descent data from $E'$ to $\Ql$
has the form
$$ \N = \langle \e, \ee, \f, \ff \rangle $$
with
$$ \NN = \langle u^s \e, u^r \ee + \e, u^s \f, u^r \ff + \f \rangle $$
and
$$ \p( u^s \e ) = b \e \,, \ \p( u^r \ee + \e ) = a \ee $$
$$ \p( u^s \f ) = b \f + v \e  \,, \ \p( u^r \ff + \f ) = a \ff + z \ee$$
with $v$, $z \in \Fl$, and generic fibre descent data satisfying
$$ [g](\e) = \gpi^{(l+1)d-ls'} \e \,, \ [g](\f)=\gpi^{(l+1)d-ls'} \f $$
$$ [g]( \ee ) = \gpi^{(l+1)c - lr'} \ee \,, \  [g]( \ff ) = \gpi^{(l+1)c 
- lr'} \ff $$
Therefore this $\Ext^1$ is at most two-dimensional over $\Fl$.
\end{thm}

\subsection{Dieudonn\'{e} module relations}  It remains to determine which
of these extensions with descent data
satisfies the relations (\ref{eqC}, \ref{eqD}, \ref{eqE})
on their Dieudonn\'{e}
module.  We check from 
the compatibility between Breuil theory and Dieudonn\'{e} theory
described in Section \ref{bm}
that each of the above extensions of Breuil modules with descent data
 yields a 
Dieudonn\'{e} module with basis
$\mathbf{v},\mathbf{w},\mathbf{v}',\mathbf{w}'$ on which $F$ and $V$
act through the matrices
$$F = 
\begin{pmatrix}
0 & 0 & 0 & 0 \\
0 & 0 & 0 & 0 \\
-b & 0 & 0 & 0 \\
 -v & -b & 0 & 0
\end{pmatrix}
$$
and
$$V = 
\begin{pmatrix}
0 & 0 & 0 & 0 \\
0 & 0 & 0 & 0 \\
1/a & 0 & 0 & 0 \\
-z/{a^2} & 1/a & 0 & 0
\end{pmatrix} .
$$ 
(Note that these matrices only describe the actions of
$F$, $V$ on this particular basis: the actions of $F$, $V$ are extended
to the full Dieudonn\'{e} module \emph{semilinearly}.)

To see this, we will have
$\mathbf{v},\mathbf{w},\mathbf{v}',\mathbf{w}'$ correspond 
respectively to the images 
of $\e,\f,\ee,\ff$ in $\N/u\N$.  Observe that $\phi(\e) = \p(u^{l^2-1} \e
) = u^{r} b\e$ which maps to $0$ in $\N/u\N$, and similarly $\phi(\f)=0$
in $\N/u\N$.  This gives the first two rows of the matrix for $F$.  Next,
$\phi(\ee) = \phi_1 (u^{l^2-1} \ee) = \phi_1 ( u^{s} (u^r \ee + \e) -
u^{s} \e) = - b \e$ in $\N/u\N$, while similarly $\phi(\ff)$ in $\N/u\N$ is  $-\phi_1(
u^{s} \f) = -b\e - v\f$.

To obtain the matrix for $V$, we note that $\p^{-1}(\e) = b^{-1}
u^{s} \e$ is $0$ in $\N/u\N$, and similar for $\p^{-1}(\f)$.  On the other
hand $\p^{-1}(\ee) = a^{-1} (u^{r} \ee + \e)$, which is $a^{-1} \e$ in
$\N/u\N$, and 
$$\p^{-1}(\ff) = a^{-1} (u^{r} \ff + \f) - \frac{z}{a^2} (u^{r} \ee +
\e)$$
which indeed is $a^{-1} \f  - \frac{z}{a^2} \e$ in $\N/u\N$. 

We know that in this case $T = {\rm Teich}(\det(\rhobar))(s)$
reduces in $\Fl$ to $ab$, and so $F + TV = 0$ precisely when $$ - v -
\frac{b}{a} z = 0 \,. $$ 

The space of extensions of Breuil modules with descent data satisfying 
the necessary 
Dieudonn\'e module relations is 
therefore at most $1$-dimensional.
This completes the proof of Theorem \ref{main}.

\nocite{Breuil}
\nocite{FontaineMazur}
\nocite{TaylorFM}
\nocite{Diamond}
\nocite{Serre}
\nocite{AST223}
\nocite{FontaineIllusie}
\nocite{Conrad}
\nocite{Mazur}
\nocite{Tate2}

\bibliographystyle{amsalpha}
\bibliography{th}

\footnotesize

\small

\vskip 0.2cm

\noindent \textsc{Department of Mathematics, McGill University, and CICMA}

\noindent \textsf{dsavitt@math.mcgill.ca}

\end{document}